\documentclass[3p]{elsarticle}
\title{The complete set of minimal simple graphs that support unsatisfiable~\mbox{\(2\)-CNFs}}
\author[uiuc]{Vaibhav Karve}
    \ead{vkarve2@illinois.edu}
\author[uiuc]{Anil N. Hirani\corref{cor1}}
    \ead{hirani@illinois.edu}

\cortext[cor1]{Corresponding author}
\address[uiuc]{Department of Mathematics, University of Illinois at
Urbana-Champaign,\\1409 W. Green Street, Urbana, IL~61801}

\usepackage{amsmath} 
\usepackage{amsfonts}
\usepackage{mathtools} 
\usepackage{color}   
\usepackage[utf8]{inputenc}
\usepackage{pstricks-add}
\usepackage{multicol} 
\usepackage{setspace}  
\usepackage[hidelinks]{hyperref}

\usepackage{tikz}
\usetikzlibrary{calc}
\usetikzlibrary{arrows.meta}
\usetikzlibrary{backgrounds}
\usetikzlibrary{decorations.pathmorphing}

\newtheorem{thm}{Theorem}
\newtheorem{lem}[thm]{Lemma}
\newtheorem{cor}[thm]{Corollary}
\newtheorem{rmk}[thm]{Remark}
\newproof{pf}{Proof}

\DeclarePairedDelimiter{\paren}{(}{)}
\DeclarePairedDelimiter\abs{\lvert}{\rvert}

\DeclarePairedDelimiterX{\Set}[2]\{\}{\, #1 \colon\, #2 \,}

\begin{document}
\newcommand{\MHG}{\mathcal{MHG}}
\newcommand{\MG}{\mathcal{MG}}
\newcommand{\HG}{\mathcal{HG}}
\newcommand{\G}{\mathcal{G}}
\newcommand{\U}{\mathcal{U}}
\newcommand{\m}{\mathrm{m}}

\newcommand{\todo}[1]{\vspace{5 mm}\par \noindent
\marginpar{\textsc{ToDo}}
\framebox{\begin{minipage}[c]{0.95 \columnwidth}
\tt #1 \end{minipage}}\vspace{5 mm}\par}

\tikzstyle{vertex}=[circle, draw, fill=black!100,
inner sep=1pt, minimum width=2pt]

\tikzstyle{vertexsmall}=[circle, draw, fill=black!100,
inner sep=0pt, minimum width=1.5pt]

\newcommand{\Pconfig}{\begin{tikzpicture}[scale=0.75, baseline=4.5ex]
\coordinate (A) at (0,0);
\coordinate (B) at (0,1);
\coordinate (C) at (30:1);
\coordinate (D) at ($(30:1) + (1,0)$);
\coordinate (E) at ($(30:2) + (1,0)$);
\coordinate (F) at ($(30:2) + (0,-1) + (1,0)$);
\draw[fill = yellow!50!white] (A) -- (B) -- (C) -- cycle;
\draw[decorate,decoration={segment length = 2mm}] (C) -- (D);
\draw[fill = yellow!50!white] (D) -- (E) -- (F) -- cycle;
\foreach \x in {(A), (B), (C), (D), (E), (F)}
\node[vertex] (\x) at \x {};
\end{tikzpicture}}

\newcommand{\Pconfigsmall}{\begin{tikzpicture}[scale=0.35, baseline=0ex]
  \coordinate (A) at (0,0);
  \coordinate (B) at (0,1);
  \coordinate (C) at (30:1);
  \coordinate (D) at ($(30:1) + (1.5,0)$);
  \coordinate (E) at ($(30:2) + (1.5,0)$);
  \coordinate (F) at ($(30:2) + (0,-1) + (1.5,0)$);
  \draw[fill = yellow!50!white] (A) -- (B) -- (C) -- cycle;
  \draw[decorate,decoration={ segment length = 1mm, amplitude=0.5mm}] (C) -- (D);
  \draw[fill = yellow!50!white] (D) -- (E) -- (F) -- cycle;
  \foreach \x in {(A), (B), (C), (D), (E), (F)}
  \node[vertexsmall] (\x) at \x {};
  \end{tikzpicture}}

\newcommand{\Vconfig}{\begin{tikzpicture}[scale=0.75, baseline=4.5ex]
\coordinate (A) at (0,0);
\coordinate (B) at (0,1);
\coordinate (C) at (30:1);
\coordinate (D) at (30:2);
\coordinate (E) at ($(30:2)+(0,-1)$);
\draw[fill = yellow!50!white] (A) -- (B) -- (C) -- cycle;
\draw[fill = yellow!50!white] (C) -- (D) -- (E) -- cycle;
\foreach \x in {(A), (B), (C), (D), (E)}
\node[vertex] (\x) at \x {};
\end{tikzpicture}}

\newcommand{\Econfig}{\begin{tikzpicture}[scale=0.75, baseline=4.5ex]
\coordinate (A) at (0,0);
\coordinate (B) at (150:1);
\coordinate (C) at (90:1);
\coordinate (D) at (30:1);
\draw[fill = yellow!50!white] (A) -- (B) -- (C) -- cycle;
\draw[fill = yellow!50!white] (A) -- (C) -- (D) -- cycle;
\foreach \x in {(A), (B), (C), (D)}
\node[vertex] (\x) at \x {};
\end{tikzpicture}}

\newcommand{\snakypath}{\begin{tikzpicture}[scale=0.75]
\coordinate (C) at (30:1);
\coordinate (D) at ($(30:1) + (1,0)$);
\draw[decorate,decoration={ segment length = 2mm}] (C) -- (D);
\end{tikzpicture}}

\tikzset{triangle/.pic={
    \coordinate (A) at (0,0);
    \coordinate (B) at (150:1);
    \coordinate (C) at (210:1);
    \draw[fill=yellow!50!white] (A)--(B)--(C)--cycle;
    \foreach \x in {(A), (B), (C)}
    \node[vertex] (\x) at \x {};
}}


\newcommand{\PPPoneconfig}{%
\begin{tikzpicture}[scale=0.5, every pic/.style={scale=0.5}, baseline=0.5ex]
\pic at (0,0) {triangle};
\pic[rotate=180] at (2,0) {triangle};
\pic[rotate=90] at (1,-1) {triangle};
\draw[decorate, decoration={segment length=1mm,amplitude=0.4mm}]
(0,0) -- (2,0) (1,-1) -- (0,0);
\end{tikzpicture}}

\newcommand{\PPPtwoconfig}{
  \begin{tikzpicture}[trim left = -12, scale=0.5, every pic/.style={scale=0.5}, baseline=0.5ex]
\pic at (0,0) {triangle};
\pic[rotate=180] at (2,0) {triangle};
\pic[rotate=90] at (1,-1) {triangle};
\draw[decorate, decoration={segment length=1mm,amplitude=0.4mm}]
(0,0) -- (2,0) (1,-1) -- (-150:1);
\end{tikzpicture}}

\newcommand{\PPVoneconfig}{
\begin{tikzpicture}[trim left = -13, scale=0.5, every pic/.style={scale=0.5}, baseline=0.5ex]
\pic at (0,0) {triangle};
\pic[rotate=180] at (0,0) {triangle};
\pic[rotate=90] at (0,-1) {triangle};
\draw[decorate, decoration={segment length=1mm,amplitude=0.4mm}]
(0,-1) -- (-150:1);
\end{tikzpicture}}

\newcommand{\PPVtwoconfig}{
\begin{tikzpicture}[trim left = -13, scale=0.5, every pic/.style={scale=0.5}, baseline=0.5ex]
\pic at (0,0) {triangle};
\pic[rotate=180] at (0,0) {triangle};
\pic[rotate=90] at (0,-1) {triangle};
\draw[decorate, decoration={segment length=1mm,amplitude=0.4mm}]
(0,0) -- (0,-1);
\end{tikzpicture}}

\newcommand{\PPEoneconfig}{
\begin{tikzpicture}[trim left = -13, scale=0.5, every pic/.style={scale=0.5}, baseline=0.5ex]
\pic[rotate=180] at ($(-1.732/2,0)$) {triangle};
\pic at ($(1.732/2,0)$) {triangle};
\pic[rotate=90] at (0,-1) {triangle};
\draw[decorate, decoration={segment length=1mm,amplitude=0.4mm}]
($(-1.732/2,0)$) -- ($(0,-1)+(-120:1)$);
\end{tikzpicture}}

\newcommand{\PPEtwoconfig}{
\begin{tikzpicture}[trim left = -13, scale=0.5, every pic/.style={scale=0.5}, baseline=0.5ex]
\pic[rotate=180] at ($(-1.732/2,0)$) {triangle};
\pic at ($(1.732/2,0)$) {triangle};
\pic[rotate=90] at (1,-1) {triangle};
\draw[decorate, decoration={segment length=1mm,amplitude=0.4mm}]
($(0,-1.732/4)$) -- ($(1,-1)$);
\end{tikzpicture}}

\newcommand{\PVVconfig}{
\begin{tikzpicture}[trim left = -17, trim right= 11, scale=0.5, every pic/.style={scale=0.5}, baseline=0.5ex]
\pic at (0,0) {triangle};
\pic[rotate=180] at (0,0) {triangle};
\pic[rotate=90] at (-150:1) {triangle};
\end{tikzpicture}}

\newcommand{\PVEconfig}{
\begin{tikzpicture}[trim left = -20, scale=0.5, every pic/.style={scale=0.5}, baseline=0.5ex]
\pic[rotate=180] at ($(-1.732/2,0)$) {triangle};
\pic at ($(1.732/2,0)$) {triangle};
\pic[rotate=90] at ($(-1.732/2,0)$) {triangle};
\end{tikzpicture}}

\newcommand{\VVVoneconfig}{
\begin{tikzpicture}[trim left = -13, scale=0.5, every pic/.style={scale=0.5}, baseline=-1.5ex]
\pic[rotate=30] at (0,0) {triangle};
\pic[rotate=150] at (0,0) {triangle};
\pic[rotate=30] at (-60:1) {triangle};
\end{tikzpicture}}

\newcommand{\VVVtwoconfig}{
\begin{tikzpicture}[trim left = -13, scale=0.5, every pic/.style={scale=0.5}, baseline=1.5ex]
\pic[rotate=-30] at (0,0) {triangle};
\pic[rotate=210] at (0,0) {triangle};
\pic[rotate=90] at (0,0) {triangle};
\end{tikzpicture}}

\newcommand{\VVEoneconfig}{
\begin{tikzpicture}[trim left = -13, scale=0.5, every pic/.style={scale=0.5}, baseline=0.5ex]
\pic[rotate=180] at ($(-1.732/2,0)$) {triangle};
\pic at ($(1.732/2,0)$) {triangle};
\begin{scope}[rotate=30,scale=1.732,transform canvas={xshift=0.433cm}]  
\coordinate (A) at (0,0);
\coordinate (B) at (150:1);
\coordinate (C) at ($(210:1)+(60:0.2)$);
\draw[fill=yellow!50!white, fill opacity=0.4, draw opacity=0.6] (A)--(B)--(C)--cycle;
\foreach \x in {(A), (B), (C)}
\node[vertex] (\x) at \x {};
\end{scope}
\end{tikzpicture}}

\newcommand{\VVEtwoconfig}{
\begin{tikzpicture}[trim left = -13, scale=0.5, every pic/.style={scale=0.5}, baseline=0.5ex]
\pic[rotate=180] at ($(-1.732/2,0)$) {triangle};
\pic at ($(1.732/2,0)$) {triangle};
\pic[rotate=90] at (0,-0.5) {triangle};:
\end{tikzpicture}}

\newcommand{\VEEconfig}{
\begin{tikzpicture}[trim left = -12.5, scale=0.5, every pic/.style={scale=0.5}, baseline=0.5ex]
\pic[rotate=180] at ($(-1.732/2,0)$) {triangle};
\pic at ($(1.732/2,0)$) {triangle};
\pic[rotate=120] at ($(-1.732/2,0)$) {triangle};
\end{tikzpicture}}

\newcommand{\EEEoneconfig}{
\begin{tikzpicture}[scale=1, baseline=4.5ex]
\coordinate (A) at (0,0);
\coordinate (B) at (1,0);
\coordinate (C) at (60:1);
\coordinate (D) at ($(C)-(0,0.52)$);
\draw[fill=yellow!50!white] (A) -- (B) -- (C) -- cycle;
\draw (A) -- (D) -- (C);
\draw (B) -- (D);
\foreach \x in {(A), (B), (C), (D)}
\node[vertex] (\x) at \x {};
\end{tikzpicture}}

\newcommand{\EEEtwoconfig}{
\begin{tikzpicture}[scale=0.57, baseline=4.5ex]
\coordinate (A) at (0,0);
\coordinate (C) at (25:1);
\coordinate (D) at ($(C)+(0,0.5)$);
\coordinate (E) at ($(C)+(0,1)$);
\coordinate (B) at ($(30:2)+(0,-1)$);
\draw[fill=yellow!50!white, fill opacity=0.3] (A) -- (B) -- (E) -- cycle;
\draw[fill=yellow!50!white, fill opacity=0.6] (A) -- (D) -- (B) -- cycle;
\draw[fill=yellow!50!white, fill opacity=1] (A) -- (C) -- (B) -- cycle;
\foreach \x in {(A), (B), (C), (D), (E)}
\node[vertex] (\x) at \x {};
\end{tikzpicture}}

\newcommand{\Kfour}{
  \begin{tikzpicture}[trim left = 2, scale=0.45, baseline=0ex]
  \coordinate (A) at (0,0);
  \coordinate (B) at (1,0);
  \coordinate (C) at (60:1);
  \coordinate (D) at ($(C)-(0,0.52)$);
  \draw[fill=yellow!50!white] (A) -- (B) -- (C) -- cycle;
  \draw (A) -- (D) -- (C);
  \draw (B) -- (D);
  \foreach \x in {(A), (B), (C), (D)}
  \node[vertexsmall] (\x) at \x {};
  \end{tikzpicture}\;}

\newcommand{\Book}{
  \begin{tikzpicture}[trim left = 2, scale=0.3, baseline=0ex]
  \coordinate (A) at (0,0);
  \coordinate (C) at (25:1);
  \coordinate (D) at ($(C)+(0,0.5)$);
  \coordinate (E) at ($(C)+(0,1)$);
  \coordinate (B) at ($(30:2)+(0,-1)$);
  \draw[fill=yellow!50!white, fill opacity=0.3] (A) -- (B) -- (E) -- cycle;
  \draw[fill=yellow!50!white, fill opacity=0.6] (A) -- (D) -- (B) -- cycle;
  \draw[fill=yellow!50!white, fill opacity=1] (A) -- (C) -- (B) -- cycle;
  \foreach \x in {(A), (B), (C), (D), (E)}
  \node[vertexsmall] (\x) at \x {};
  \end{tikzpicture}\;}

\newcommand{\Butterfly}{
\begin{tikzpicture}[trim left = 1.5, scale=0.35, baseline=0ex]
  \coordinate (A) at (0,0);
  \coordinate (B) at (0,1);
  \coordinate (C) at (30:1);
  \coordinate (D) at (30:2);
  \coordinate (E) at ($(30:2)+(0,-1)$);
  \draw[fill=yellow!50!white] (A) -- (B) -- (C) -- cycle;
  \draw[fill=yellow!50!white] (C) -- (D) -- (E) -- cycle;
  \foreach \x in {(A), (B), (C), (D), (E)}
  \node[vertexsmall] (\x) at \x {};
  \end{tikzpicture}\;}

\newcommand{\Butterflybig}{
  \begin{tikzpicture}[scale=0.7, baseline=3.5ex]
  \coordinate (A) at (0,0);
  \coordinate (B) at (0,1);
  \coordinate (C) at (30:1);
  \coordinate (D) at ($(30:1)$);
  \coordinate (E) at ($(30:2)$);
  \coordinate (F) at ($(30:2) + (0,-1)$);
  \draw[fill = yellow!50!white] (A) -- (B) -- (C) -- cycle;
  \draw[fill = yellow!50!white] (D) -- (E) -- (F) -- cycle;
  \foreach \x in {(A), (B), (C), (D), (E), (F)}
  \node[vertex] (\x) at \x {};
  \end{tikzpicture}}

\newcommand{\Bowtie}{
  \begin{tikzpicture}[trim left = 1.5, scale=1, baseline=0ex]
  \coordinate (A) at (0,0);
  \coordinate (B) at (0,1);
  \coordinate (C) at (30:1);
  \coordinate (D) at ($(30:1) + (1,0)$);
  \coordinate (E) at ($(30:2) + (1,0)$);
  \coordinate (F) at ($(30:2) + (1,-1)$);
  \draw[fill=yellow!00!white] (A) -- (B) -- (C) -- cycle;
  \draw[fill=yellow!00!white] (C) -- (D);
  \draw[fill=yellow!00!white] (D) -- (E) -- (F) -- cycle;
  \foreach \x in {(A), (B), (C), (D), (E), (F)}
  \node[vertex] (\x) at \x {};
  \end{tikzpicture}\;}

\newcommand{\Bowtiesmall}{
  \begin{tikzpicture}[trim left = 1.5, scale=0.35, baseline=0ex]
  \coordinate (A) at (0,0);
  \coordinate (B) at (0,1);
  \coordinate (C) at (30:1);
  \coordinate (D) at ($(30:1) + (1,0)$);
  \coordinate (E) at ($(30:2) + (1,0)$);
  \coordinate (F) at ($(30:2) + (1,-1)$);
  \draw[fill=yellow!50!white] (A) -- (B) -- (C) -- cycle;
  \draw[fill=yellow!00!white] (C) -- (D);
  \draw[fill=yellow!50!white] (D) -- (E) -- (F) -- cycle;
  \foreach \x in {(A), (B), (C), (D), (E), (F)}
  \node[vertexsmall] (\x) at \x {};
  \end{tikzpicture}\;}


\newcommand{\KfourMinusE}{
  \begin{tikzpicture}[trim left = 2, trim right =9, scale=0.35, baseline=0ex]
  \coordinate (A) at (0,0);
  \coordinate (B) at (1,0);
  \coordinate (C) at (1,1);
  \coordinate (D) at (0,1);
  \draw[fill=yellow!50!white] (A) -- (B) -- (C) -- (D) -- cycle;
  \draw (A) -- (C);
  \foreach \x in {(A), (B), (C), (D)}
  \node[vertexsmall] (\x) at \x {};
  \end{tikzpicture}\;}

\newcommand{\BookSubgraphOne}{
  \begin{tikzpicture}[trim left = 2, scale=0.3, baseline=0ex]
  \coordinate (A) at (0,0);
  \coordinate (C) at (25:1);
  \coordinate (D) at ($(C)+(0,0.5)$);
  \coordinate (E) at ($(C)+(0,1)$);
  \coordinate (B) at ($(30:2)+(0,-1)$);
  \draw (B) -- (E);
  \draw[fill=yellow!50!white, fill opacity=0.6] (A) -- (D) -- (B) -- cycle;
  \draw[fill=yellow!50!white, fill opacity=1]   (A) --(C) -- (B) -- cycle;
  \foreach \x in {(A), (B), (C), (D), (E)}
  \node[vertexsmall] (\x) at \x {};
  \end{tikzpicture}\;}

\newcommand{\BookSubgraphTwo}{
  \begin{tikzpicture}[trim left = 2, scale=0.3, baseline=0ex]
  \coordinate (A) at (0,0);
  \coordinate (C) at (25:1);
  \coordinate (D) at ($(C)+(0,0.5)$);
  \coordinate (E) at ($(C)+(0,1)$);
  \coordinate (B) at ($(30:2)+(0,-1)$);
  \draw[fill=yellow!00!white] (A) -- (E) -- (B)  -- (C) -- cycle;
  \draw (A) -- (E) -- (B);
  \draw (A) -- (D) -- (B);
  \draw (A) -- (C) -- (B);
  \foreach \x in {(A), (B), (C), (D), (E)}
  \node[vertexsmall] (\x) at \x {};
  \end{tikzpicture}\;}

\newcommand{\ButterflySubgraphOne}{
  \begin{tikzpicture}[trim left = 0, scale=0.35, baseline=0ex]
  \coordinate (A) at (0,0);
  \coordinate (B) at (0,1);
  \coordinate (C) at (30:1);
  \coordinate (D) at (30:2);
  \coordinate (E) at ($(30:2)+(0,-1)$);
  \draw[fill=yellow!50!white] (A) -- (B) -- (C) -- cycle;
  \draw[fill=yellow!00!white] (C) -- (D);
  \draw (C) -- (E);
  \foreach \x in {(A), (B), (C), (D), (E)}
  \node[vertexsmall] (\x) at \x {};
  \end{tikzpicture}\;}

\newcommand{\ButterflySubgraphTwo}{
  \begin{tikzpicture}[trim left = 0, scale=0.35, baseline=0ex]
  \coordinate (A) at (0,0);
  \coordinate (B) at (0,1);
  \coordinate (C) at (30:1);
  \coordinate (D) at (30:2);
  \coordinate (E) at ($(30:2)+(0,-1)$);
  \draw[fill=yellow!50!white] (A) -- (B) -- (C) -- cycle;
  \draw[fill=yellow!00!white] (C) -- (D) -- (E);
  \foreach \x in {(A), (B), (C), (D), (E)}
  \node[vertexsmall] (\x) at \x {};
  \end{tikzpicture}\;}

\newcommand{\BowtieSubgraphOne}{
  \begin{tikzpicture}[trim left = 0, scale=0.35, baseline=0ex]
  \coordinate (A) at (0,0);
  \coordinate (B) at (0,1);
  \coordinate (C) at (30:1);
  \coordinate (D) at ($(30:1)+(1,0)$);
  \coordinate (E) at ($(30:2)+(1,0)$);
  \coordinate (F) at ($(30:2)+(1,0)+(0,-1)$);
  \draw[fill=yellow!50!white] (A) -- (B) -- (C) -- cycle;
  \draw[fill=yellow!00!white] (D) -- (E) -- (F);
  \draw (C) -- (D);
  \foreach \x in {(A), (B), (C), (D), (E), (F)}
  \node[vertexsmall] (\x) at \x {};
  \end{tikzpicture}\;}

\newcommand{\BowtieSubgraphTwo}{
  \begin{tikzpicture}[trim left = 0, scale=0.35, baseline=0ex]
  \coordinate (A) at (0,0);
  \coordinate (B) at (0,1);
  \coordinate (C) at (30:1);
  \coordinate (D) at ($(30:1)+(1,0)$);
  \coordinate (E) at ($(30:2)+(1,0)$);
  \coordinate (F) at ($(30:2)+(1,0)+(0,-1)$);
  \draw[fill=yellow!50!white] (A) -- (B) -- (C) -- cycle;
  \draw[fill=yellow!00!white] (D) -- (E);
  \draw (C) -- (D);
  \draw (D) -- (F);
  \foreach \x in {(A), (B), (C), (D), (E), (F)}
  \node[vertexsmall] (\x) at \x {};
  \end{tikzpicture}\;}

\newcommand{\BowtieSubgraphThree}{
  \begin{tikzpicture}[trim left = 0, scale=0.35, baseline=0ex]
  \coordinate (A) at (0,0);
  \coordinate (B) at (0,1);
  \coordinate (C) at (30:1);
  \coordinate (D) at ($(30:1)+(1,0)$);
  \coordinate (E) at ($(30:2)+(1,0)$);
  \coordinate (F) at ($(30:2)+(1,0)+(0,-1)$);
  \draw[fill=yellow!50!white] (A) -- (B) -- (C) -- cycle;
  \draw[fill=yellow!50!white] (D) -- (E) -- (F) -- cycle;
  \foreach \x in {(A), (B), (C), (D), (E), (F)}
  \node[vertexsmall] (\x) at \x {};
  \end{tikzpicture}\;}

\newcommand{\FlexibleButterfly}{
  \begin{tikzpicture}[scale=1, baseline=0ex]
  \coordinate (A) at (0,0);
  \coordinate (TL) at (-1,0.5);
  \coordinate (BR) at (1,-0.5);
  \clip (TL) rectangle (BR);
  \draw (A) to [out=120, in=240, loop, min distance=20mm] (A);
  \draw (A) to [out=60, in=300, loop, min distance=20mm] (A);
  \foreach \x in {(A)}
  \node[vertex] (\x) at \x {};
  \end{tikzpicture}\;}

\newcommand{\FlexibleBowtie}{
  \begin{tikzpicture}[scale=1, baseline=0ex]
  \coordinate (A) at (0,0);
  \coordinate (B) at (0.75,0);
  \coordinate (TL) at (-1.0,0.5);
  \coordinate (BR) at (1.75,-0.5);
  \clip (TL) rectangle (BR);
  \draw (A) to [out=120, in=240, loop, min distance=20mm] (A);
  \draw (B) to [out=60, in=300, loop, min distance=20mm] (B);
  \draw (A) -- (B);
  \foreach \x in {(A), (B)}
  \node[vertex] (\x) at \x {};
  \end{tikzpicture}\;}

\newcommand{\FlexibleKfour}{
  \begin{tikzpicture}[scale=1, baseline=0ex]
  \coordinate (A) at (0,0);
  \coordinate (B) at (0,0.5);
  \coordinate (C) at ($(0,0) + (-30:0.5)$);
  \coordinate (D) at ($(0,0) + (-150:0.5)$);
  \coordinate (TL) at (-1,0.7);
  \coordinate (BR) at (1,-0.9);
  \clip (TL) rectangle (BR);
  \draw (A) -- (B);
  \draw (A) -- (C);
  \draw (A) -- (D);
  \draw (B) to [out=30, in=30, min distance=7mm] (C);
  \draw (B) to [out=150, in=150, min distance=7mm] (D);
  \draw (D) to [out=270, in=270, min distance=7mm] (C);
  
  \foreach \x in {(A), (B), (C), (D)}
  \node[vertex] (\x) at \x {};
  \end{tikzpicture}\;}

\newcommand{\FlexibleBook}{
  \begin{tikzpicture}[scale=1, baseline=3.5ex]
  \coordinate (A) at (0,0);
  \coordinate (B) at (1.25,0);
  \coordinate (TL) at (-0.25,1.25);
  \coordinate (BR) at (1.5,-0.35);
  \clip (TL) rectangle (BR);
  \draw (A) to [out=90, in=90, min distance=15mm] (B);
  \draw (A) to [out=70, in=110, min distance=12mm] (B);
  \draw (A) to [out=50, in=130, min distance=9mm] (B);
  \draw (A) -- (B);
  \foreach \x in {(A), (B)}
  \node[vertex] (\x) at \x {};
  \end{tikzpicture}\;}

\newcommand{\WhiteButterfly}{
  \begin{tikzpicture}[trim left = 1.5, scale=0.35, baseline=0ex]
  \coordinate (A) at (0,0);
  \coordinate (B) at (0,1);
  \coordinate (C) at (30:1);
  \coordinate (D) at (30:2);
  \coordinate (E) at ($(30:2)+(0,-1)$);
  \draw[fill=yellow!00!white] (A) -- (B) -- (C) -- cycle;
  \draw[fill=yellow!00!white] (C) -- (D) -- (E) -- cycle;
  \foreach \x in {(A), (B), (C), (D), (E)}
  \node[vertexsmall] (\x) at \x {};
  \end{tikzpicture}\;}

\newcommand{\SquareButterfly}{
  \begin{tikzpicture}[trim left = 1.5, scale=0.35, baseline=0ex]
  \coordinate (A) at (0,0);
  \coordinate (B) at (0,1);
  \coordinate (C) at (1,0);
  \coordinate (D) at (1,1);
  \coordinate (E) at (2,0);
  \coordinate (F) at (2,1);
  \draw[fill=yellow!00!white] (A) -- (B) -- (D) -- (C) -- cycle;
  \draw[fill=yellow!00!white] (C) -- (D) -- (F) -- (E) -- cycle;
  \foreach \x in {(A), (B), (C), (D), (E), (F)}
  \node[vertexsmall] (\x) at \x {};
  \end{tikzpicture}\;}

\newcommand{\Twohills}{
  \begin{tikzpicture}[scale=0.8, baseline=1.5ex]
  \coordinate (A) at (0,0);
  \coordinate (B) at (1,0);
  \coordinate (C) at (60:1);
  \coordinate (D) at (2,0);
  \coordinate (E) at ($(1,0) + (60:1)$);
  \draw[fill = yellow!00!white] (A) -- (B) -- (C) -- cycle;
  \draw[fill = yellow!00!white] (B) -- (D) -- (E) -- cycle;
  \foreach \x in {(A), (B), (C), (D), (E)}
  \node[vertex] (\x) at \x {};
  \end{tikzpicture}}

\newcommand{\Threehills}{
  \begin{tikzpicture}[scale=0.8, baseline=1.5ex]
  \coordinate (A) at (0,0);
  \coordinate (B) at (1,0);
  \coordinate (C) at (60:1);
  \coordinate (D) at (2,0);
  \coordinate (E) at ($(1,0) + (60:1)$);
  \coordinate (F) at (3,0);
  \coordinate (G) at ($(2,0) + (60:1)$);
  \draw[fill = yellow!00!white] (A) -- (B) -- (C) -- cycle;
  \draw[fill = yellow!00!white] (B) -- (D) -- (E) -- cycle;
  \draw[fill = yellow!00!white] (D) -- (F) -- (G) -- cycle;
  \foreach \x in {(A), (B), (C), (D), (E), (F), (G)}
  \node[vertex] (\x) at \x {};
  \end{tikzpicture}}

\newcommand{\Fourhills}{
  \begin{tikzpicture}[scale=0.8, baseline=1.5ex]
  \coordinate (A) at (0,0);
  \coordinate (B) at (1,0);
  \coordinate (C) at (60:1);
  \coordinate (D) at (2,0);
  \coordinate (E) at ($(1,0) + (60:1)$);
  \coordinate (F) at (3,0);
  \coordinate (G) at ($(2,0) + (60:1)$);
  \coordinate (H) at (4,0);
  \coordinate (I) at ($(3,0) + (60:1)$);
  \draw[fill = yellow!00!white] (A) -- (B) -- (C) -- cycle;
  \draw[fill = yellow!00!white] (B) -- (D) -- (E) -- cycle;
  \draw[fill = yellow!00!white] (D) -- (F) -- (G) -- cycle;
  \draw[fill = yellow!00!white] (F) -- (H) -- (I) -- cycle;
  \foreach \x in {(A), (B), (C), (D), (E), (F), (G), (H), (I)}
  \node[vertex] (\x) at \x {};
  \end{tikzpicture}}

\begin{abstract}
  A propositional logic sentence in conjunctive normal form that has
  clauses of length two (a 2-CNF) can be associated with a
  multigraph in which the vertices correspond to the variables and
  edges to clauses. We first show that every such sentence that has
  been reduced, that is, which is unchanged under application of
  certain tautologies, is equisatisfiable to a 2-CNF whose associated
  multigraph is, in fact, a simple graph. Our main result is a
  complete characterization of graphs that can support unsatisfiable
  2-CNF sentences. We show that a simple graph can support an
  unsatisfiable reduced 2-CNF sentence if and only if it contains any
  one of four specific small graphs as a topological minor.
  Equivalently, all reduced 2-CNF sentences supported on a
  given simple graph are satisfiable if and only if all subdivisions
  of those four graphs are forbidden as subgraphs of of the original
  graph.
  We conclude with a discussion of why the Robertson-Seymour graph
  minor theorem does not apply in our approach.
\end{abstract}
\maketitle

\section{Introduction}\label{sec:intro}
Given a sentence in propositional logic, the satisfiability decision
problem is to determine if there exists a truth assignment for the
variables that makes the sentence true. Let~\(V\) be a finite set of
Boolean variables,~\(\neg V\) be the set~\(\Set*{\neg x}{x\in V}\) of
negations, and let the symbols~\(\top\) and~\(\bot\) be \emph{True}
and \emph{False} values of variables.  We define the \emph{set of
  literals obtained from~\(V\)} to be the
set~\(V\cup\neg V\cup\{\top,\bot\}\).

A \emph{Conjunctive Normal Form (CNF) on~\(V\)} is a conjunction of
one or more clauses, where each \emph{clause} is a disjunction of
literals.  A CNF (or clause) is \emph{reduced} if it is unchanged
under application of all the tautologies listed below
\begin{multicols}{4}
\begin{enumerate}[]
  \item \(x\vee\top = \top\),
  \item \(x\vee\bot = x\),
  \item \(x\wedge\top = x\),
  \item \(x\wedge\bot = \bot\),
  \item \(x\vee x = x\),
  \item \(x\wedge x = x\),
  \item \(\neg(\neg x) = x\),
  \item \(x\vee\neg x = \top\).
\end{enumerate}
\end{multicols}
It is a fact that every CNF (or clause) is logically equivalent to a
reduced CNF (or clause).  The \emph{length} of a reduced clause is the
number of literals in the clause.  For~\(k\in\mathbb{N}\), a reduced
CNF is a \emph{reduced~\(k\)-CNF} if all of its clauses have
length~\(k\).  Additionally, for ease of stating results, the
CNFs~\(\top\) and~\(\bot\) will be treated as~\mbox{\(k\)-CNFs} for
every~\(k\in\mathbb{N}\).  We will refer to the satisfiability
decision problem as SAT and the satisfiability problem for
$k$-CNFs as $k$-SAT.

Instead of algorithmic issues our aim in this paper is to study the
structure of unsatisfiable 2-CNF sentences in a sense that will be
made precise later. For completeness, here we briefly summarize the
relevant fundamental algorithmic results for satisfiability
problems. The 2-SAT problem is in P. An $\mathcal{O}(n^4)$
algorithm was given by Krom~\cite{Krom1967} and a \mbox{linear-time} algorithm by Even,
Itai and Shamir~\cite{EvItSh1976} and Aspvall, Plass and
Tarjan~\cite{AsPlTa1979}. All solutions of a given 2-CNF sentence can
be listed efficiently using an algorithm by Feder~\cite{Feder1994}.

In contrast with the algorithmic tractability of 2-SAT, the SAT
problem in general is NP-complete as was shown by Cook and by Levin
independently~\cite{Cook1971, Levin1973}. As part of the proof of the
NP-completeness of SAT, they also proved that every logical sentence
can be rewritten as a CNF while changing its length by no more than a
constant factor. Schaefer's dichotomy theorem states necessary and
sufficient conditions under which a finite set \(S\) of relations over
the Boolean domain yields \mbox{polynomial-time} or \mbox{NP-complete} problems
when the relations of \(S\) are used to constrain some of the propositional
variables~\cite{Schaefer1978}. Thus~\cite{Schaefer1978} gives a
necessary and sufficient condition for \mbox{SAT-type} problems to be in P
vs.\ NP.

Our paper relates properties of certain graphs to satisfiability. An
early connection between satisfiability and graphs was in the proof of
NP-completeness of various graph problems, such as the clique decision
problem and the vertex cover problem, by Karp \cite{Karp1972}. One
of the \mbox{linear-time} algorithms for 2-SAT~\cite{AsPlTa1979} mentioned
above also related graphs and satisfiability in its
use of \mbox{strongly-connected} graph components as a tool for deciding
satisfiability.

We explore the structures of unsatisfiable 2-CNFs by relating reduced
2-CNFs to graphs and examining which graphs can support unsatisfiable
sentences. Given a 2-CNF, an associated multigraph can be created by
identifying the variables as vertices and each clause as an
edge. Since multiple clauses may involve the same two variables, it is
not immediate that it is sufficient to consider graphs rather than
multigraphs. That it indeed is so is the content of
Section~\ref{sec:simple}. In Section~\ref{sec:homeo} we theorems about
the relation between graphs that support unsatisfiable
sentences. Theorem~\ref{thm: U_closure} in that section shows that the
family of simple graphs that can support an unsatisfiable sentence is
closed under graph homeomorphism while Theorem~\ref{thm: subgraph_U}
shows that a graph can support an unsatisfiable sentence if one of its
subgraph can. Theorem~\ref{thm: edge_contraction_U} shows that if a graph
can support an unsatisfiable sentence, then the graphs obtained by
\mbox{edge-contractions} at edges not contained in triangles can.
Section~\ref{sec:structure} is about connectivity properties of
graphs that we need to prove the main result. The main result of this
paper is Theorem~\ref{thm:main} in which we give a complete
characterization of graphs that can support unsatisfiable
sentences. This is given in the form of a finite set of obstructions
to supporting only satisfiable sentences. In Section~\ref{sec:conclusion}
we discuss how our approach differs from the application of finite
obstructions (forbidden minors) theory developed in the Robertson-Seymour
graph minor theorem published in a series of papers starting
with~\cite{RoSe1983} and ending with~\cite{RoSe2004}.

\section{Preliminaries}\label{sec:prelim}

For a reduced CNF~\(S\) and a variable~\(v\), we denote
by~\(S[v:=\top]\) the reduced CNF obtained by
\begin{enumerate}[a)]
  \item setting all occurrences of~\(v\) in~\(S\) to~\(\top\),
  \item setting all occurrences of~\(\neg v\) in~\(S\) to~\(\bot\),
    and
  \item reducing the CNF thus obtained, that is, applying all the
    tautologies listed in Section~\ref{sec:intro} to the CNF in order to
    obtain a reduced CNF.
\end{enumerate}
Similarly, the reduced CNF obtained by setting~\(v\) to~\(\bot\),
setting~\(\neg v\) to~\(\top\), and then performing the tautological
reductions is denoted by~\(S[v:=\bot]\) .  If we set multiple
variables to true or false, we use the
notation~\({S[\{a,b\}:=\top;\;\{c,d\} :=\bot]}\) as
\mbox{short-hand} for~\({S[a:=\top][b:=\top][c:=\bot][d:=\bot]}\)
and so on for variables~\({a,b,c}\) and~\(d\).  Note that the order
in which we set variables to true or false is not important, as it
yields logically equivalent results.

We now define some more adjectives for CNFs.  A reduced CNF~\(S\) is
\emph{true} if~\mbox{\(S = \top\)}, \emph{false} if~\(S = \bot\), and
\emph{nontrivial} if~\(S\) is neither true nor false.  A reduced
CNF~\(S\) is \emph{satisfiable} if there is a subset~\(V_1\) of~\(V\)
such that~\mbox{\(S[V_1:=\top;V - V_1:=\bot]\)} is true.  A reduced
CNF is \emph{unsatisfiable} if it is not satisfiable.  Two reduced
CNFs,~\(S\) and~\(S'\) are \emph{equisatisfiable}, denoted
by~\(S\sim S'\), if either both are satisfiable or both are
unsatisfiable.

\section{Graph associated with a reduced
  CNF}\label{sec:associated_graphs}

We define an operation~\({\abs*{{}\cdot{}}: V\cup\neg V\rightarrow V}\)
as~\({\abs*{x}= \abs*{\neg x}= x}\), for every~\({x\in V}\).
Let~\(\bigvee_{i=1}^k a_i\) be a nontrivial clause denoted
by~\(\alpha\), where~\(k\in\mathbb{N}\) and~\(a_i\in V\cup\neg V\)
for~\({i\in\{1,\ldots,k\}}\).  We extend the domain of the
operation~\(\abs*{{}\cdot{}}\) to include clauses like~\(\alpha\) by
defining~\({\abs*{\alpha} = \{\abs*{a_i} : i\in\{1,\ldots,k\}\}}\).

Let~\(S\) be a nontrivial~\mbox{\(2\)-CNF}.  We can write it
as~\(S = \bigwedge_{j=1}^n\alpha_j\), where \(n\in\mathbb{N}\) and
each \(\alpha_j\) is a clause of \mbox{length \(2\)}.  To \(S\) we
associate a multigraph \(\MG(S)\) having vertex set
\(\bigcup_{j=1}^n\abs*{\alpha_j}\) and having an edge
\(\abs*{\alpha_j}\) for every \(j\in\{1,\ldots,n\}\).  Note that a
variable \(x\) and its negation \(\neg x\) are both represented by a
single vertex (also labeled~\(x\)) in \(\MG(S)\).  The
\mbox{\(2\)-CNF} \(S\) being nontrivial guarantees that \(\MG(S)\) is
a nonempty multigraph, thus avoiding pathological cases.  The
nontriviality of \(S\) implies that \(S\) is reduced, thus ensuring
there are no \mbox{self-loops} in~\(\MG(S)\), since \mbox{self-loops}
can only occur due to clauses~\((x\vee \neg x)\) and~\mbox{\((x\vee x)\)}.

If a multigraph \(G\) is the associated multigraph of a
\mbox{\(2\)-CNF} \(S\), then we say~\(S\) is \emph{supported on
  \(G\)}.  A nontrivial \mbox{\(2\)-CNF} \(S\) is \emph{simple} if its
associated graph is simple.  In such a case, we drop the prefix
``multi-'' and denote the associated graph \(\MG(S)\) simply
by~\(\G(S)\). We denote by~\(\U\) the following family of graphs
\begin{equation*}
  \U = \Set*{\G(S)}{S \mbox{ is an unsatisfiable simple } 2\mbox{-CNF}}.
\end{equation*}
Our aim in this paper is to characterize the elements of~\(\U\).
First, we note that~\(\U\) is nonempty since
  \[S = (a\vee b) \wedge(\neg a\vee c)
         \wedge(\neg b\vee c) \wedge(\neg c\vee d)
         \wedge(\neg d\vee e) \wedge(\neg d\vee f)
         \wedge(\neg e\vee \neg f),\]
is an unsatisfiable simple \mbox{\(2\)-CNF}
supported on the following graph.
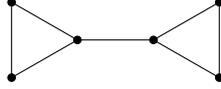
\[\Bowtie\]

\section{Simple CNFs suffice} \label{sec:simple}

In this section, we show that every nontrivial
reduced~\mbox{\(2\)-CNF} is equisatisfiable to a
simple~\mbox{\(2\)-CNF}.  Thus when studying satisfiability
of~\mbox{\(2\)-CNFs} we only need consider those that are simple.

First, we make the observation that for any~\(a\in V\), there are at
most two length~\(1\) clauses, namely~\(a\) and~\(\neg a\).  Such
clauses involving a single variable~\(a\) (or its negation) will be
referred to as~\mbox{\emph{\((a)\)-clauses}}.  For every~\(a,b\in V\),
there are at most four length~\(2\) clauses,
namely~\((a\vee b)\),~\((a\vee\neg b)\),~\((\neg a\vee b)\)
and~\((\neg a\vee\neg b)\).  Consequently, for any
reduced~\mbox{\(2\)-CNF}~\(S\), edges in~\(\MG(S)\) have a maximum
multiplicity of~\(4\).  Such clauses involving \emph{both}~\(a\) and~\(b\)
(or their negations) will be referred to
as~\mbox{\emph{\((a,b)\)-clauses}}.

\begin{lem}\label{lem: a_clause_cases}
  Let~\(S\) and~\(R\) be reduced CNFs.  Let~\(a\in V\) such that~\(R\)
  does not contain any~\mbox{\((a)\)-clauses}.
  \begin{enumerate}[a)]
  \item If~\(S = R\wedge a\wedge\neg a\), then~\(S\) is
    unsatisfiable. \label{item: a_clause_2}
  \item If~\(S = R\wedge a\), then~\(S\sim R[a:=T]\). \label{item:
      a_clause_1}
  \end{enumerate}
\end{lem}
\begin{pf}
  Let~\({S = R\wedge a\wedge\neg a}\).  We can simplify this
  as~\({S = R\wedge\bot = \bot}\).  Therefore,~\(S\) is unsatisfiable.

  Let~\({S = R\wedge a}\).  Any truth assignment that satisfies~\(S\)
  must have~\(\left[a:=\top\right]\) and must therefore
  satisfy~\(R\left[a:=\top\right]\).  Conversely, consider a truth
  assignment that satisfies~\({R\left[a:=\top\right]}\).  Recall
  that~\({R\left[a:=\top\right]}\) denotes the reduced CNF resulting
  from replacing~\(a\) with~\(\top\) and~\(\neg a\) with~\(\bot\)
  in~\(R\).  Hence, the variable~\(a\) (or its negation) are not
  present in~\({R\left[a:=\top\right]}\).  As a result, any truth
  assignment that satisfies~\(R\left[a:=\top\right]\) must leave the
  variable~\(a\) unassigned.  By additionally
  setting~\(\left[a:=\top\right]\), we obtain an assignment that
  satisfies~\(S\).  We conclude
  that~\({S\sim R\left[a:=\top\right]}\).  \qed\end{pf}

\begin{lem}\label{lem: ab_clause_cases}
  Let~\(S\) and~\(R\) be reduced~\mbox{\(2\)-CNFs}.  Let~\(a,b\in V\)
  such that~\(R\) does not contain any~\mbox{\((a,b)\)-clauses}.
  \begin{enumerate}[a)]
  \item
    If~\({S = R\wedge(a\vee b)\wedge(a\vee\neg b)\wedge(\neg a\vee
      b)\wedge(\neg a\vee\neg b)}\), then~\(S\) is
    unsatisfiable.\label{item: ab_4_clause}
  \item
    If~\({S = R\wedge(a\vee b)\wedge(a\vee \neg b)\wedge(\neg a\vee
      b)}\), then~\({S\sim R[\{a,b\}:=\top]}\). \label{item:
      ab_3_clause}
  \item If~\({S = R\wedge(a\vee b)\wedge(a\vee \neg b)}\),
    then~\({S\sim R[a:=\top]}\). \label{item: ab_2_clause_1}
  \item If~\({S = R\wedge(a\vee b)\wedge(\neg a\vee b)}\),
    then~\({S\sim R[b:=\top]}\). \label{item: ab_2_clause_2}
  \item If~\({S = R\wedge(a\vee b)\wedge(\neg a\vee \neg b)}\),
    then~\({S\sim R[b:=\neg a]}\) \label{item: ab_2_clause_3}
  \end{enumerate}
\end{lem}
\begin{pf}
  Suppose~\(S\) is of the
  form~\({R\wedge(a\vee b)\wedge(a\vee\neg b)\wedge(\neg a\vee
    b)\wedge(\neg a\vee\neg b)}\). As a result of the logical
  implications~\({(a\vee b)\wedge(a\vee\neg b)\rightarrow a}\;\) and
  \\\({(\neg a\vee b)\wedge(\neg a\vee\neg b)\rightarrow\neg a}\),
  we can deduce the
  implication~\({S\rightarrow a\wedge \neg a = \bot}\).
  Therefore,~\(S\) is unsatisfiable.

  Let~\({S = R\wedge(a\vee b)\wedge(a\vee \neg b)\wedge(\neg a\vee
    b)}\).  The three~\mbox{\((a,b)\)-clauses} of~\(S\)
  imply~\((a\wedge b)\).  Hence any truth assignment that
  satisfies~\(S\) must also satisfy~\({(a\wedge b)}\).  We can
  conclude that such an assignment must have~\([\{a,b\}:=\top]\).
  Hence this assignment must also
  satisfy~\(R\left[\{a,b\}:=\top\right]\).  Thus,~\(S\) satisfiable
  implies \(R\left[\{a,b\}:=\top\right]\) satisfiable.  Conversely,
  consider a truth assignment that
  satisfies~\({R\left[\{a,b\}:=\top\right]}\).  By additionally
  setting~\([\{a,b\}:=\top]\), we obtain an assignment that
  satisfies~\(S\).  We conclude~\({S\sim R[\{a,b\}:=\top]}\).

  Let~\({S = R\wedge(a\vee b)\wedge(a\vee \neg b)}\).
  Since~\({(a\vee b)\wedge(a\vee\neg b)\rightarrow a}\), any truth
  assignment that satisfies~\(S\) must have~\(\left[a:=\top\right]\).
  Thus, such a truth assignment must also satisfy~\(R[a:=\top]\).
  Conversely, consider a truth assignment that
  satisfies~\(R[a:=\top]\).  By additionally setting~\([a:=\top]\), we
  obtain an assignment that satisfies~\(S\).  We conclude
  that~\(S\sim R[a:=\top]\).

  Let~\({S = R\wedge(a\vee b)\wedge(\neg a\vee b)}\).
  Since~\({(a\vee b)\wedge(\neg a\vee b)\rightarrow b}\), by a similar
  reasoning to the previous case, we conclude
  that~\(S\sim R[b:=\top]\).

  Lastly, let~\({S = R\wedge(a\vee b)\wedge(\neg a\vee \neg b)}\).
  Any truth assignment that satisfies~\(S\) must map~\(a\)
  and~\(b\) to opposite truth values in order to satisfy the
  two~\mbox{\((a,b)\)-clauses}.  Thus, such a truth assignment must
  also satisfy~\(R[b:=\neg a]\).  Conversely, consider a truth
  assignment that satisfies~\(R[b:=\neg a]\).  In this truth
  assignment, we have no assignment for \(b\) (as all occurrences of
  \(b\) have been replaced by \(\neg a\)).  By additionally
  setting~\([b:=\neg a]\), we obtain an assignment that
  satisfies~\(S\).  We conclude that~\({S\sim R[b:=\neg a]}\).
  \qed\end{pf}

\begin{lem}\label{lem: reduction_of_multiplicities}
  Let~\(S\) be a reduced CNF with clauses of length at most~\(2\).
  \begin{enumerate}[a)]
  \item If~\(S\) has an~\mbox{\((a)\)-clause} for some variable~\(a\),
    then there exists a reduced \mbox{\(2\)-CNF} \(S'\)
    equisatisfiable to~\(S\).
  \item If~\(S\) has four~\mbox{\((a,b)\)-clauses} for
    some~\(a,b\in V\), then~\(S\) is unsatisfiable.
  \item If~\(S\) has three or fewer~\mbox{\((a,b)\)-clauses} for
    every~\(a,b\in V\), then there exists a~\mbox{\(2\)-CNF}~\(S'\)
    equisatisfiable to~\(S\) such that either~\(S'\) is trivial
    or~\(S'\) is simple.
  \end{enumerate}
\end{lem}
\begin{pf}
  Let~\(S\) be a reduced CNF having four~\mbox{\((a,b)\)-clauses} for
  some~\(a,b\in V\).  We can then write~\(S\) in the
  form
    \[S = R\wedge(a\vee b)\wedge(a\vee\neg b)\wedge(\neg a\vee
    b)\wedge(\neg a\vee\neg b),\]
  where~\(R\) is a reduced~\mbox{\(2\)-CNF} not containing
  any~\mbox{\((a,b)\)-clauses}.  By Lemma~\ref{lem: ab_clause_cases}~\eqref{item: ab_4_clause}, we conclude that~\(S\)
  is unsatisfiable.
  
  Let~\(S\) have an~\mbox{\((a)\)-clause} for some~\(a\in V\).  The
  reduced CNF~\(S\) has either exactly one or exactly
  two~\mbox{\((a)\)-clauses}.  If~\(S\) has exactly
  two $(a)$-clauses, then we can
  write~\({S = R\wedge a\wedge \neg a}\), where~\(R\) is a reduced CNF
  not containing any~\mbox{\((a)\)-clauses}.  By Lemma~\ref{lem:
    a_clause_cases}~\eqref{item: a_clause_2}, we conclude that~\(S\)
  is unsatisfiable.  Hence we can write~\(S\sim S_1\), where~\(S_1\)
  is the trivially false~\mbox{\(2\)-CNF}~\(\bot\).

  If~\(S\) has exactly one~\mbox{\((a)\)-clause}, then by
  exchanging~\(\neg a\) with~\(a\) if needed, we can ensure that~\(S\)
  is of the form~\(R\wedge a\), where~\(R\) is a reduced CNF not
  containing any~\((a)\)-clauses.  By Lemma~\ref{lem:
    a_clause_cases}~\eqref{item: a_clause_1}, we conclude
  that~\(S\sim R[a:=\top]\).  The reduced CNF~\(R[a:=\top]\) has
  no~\mbox{\((a)\)-clauses} but might have either one or
  two~\mbox{\((b)\)-clauses} for some other variable~\(b\).  By
  repeated application of Lemma~\ref{lem:
    a_clause_cases}~\eqref{item: a_clause_2} and~\eqref{item:
    a_clause_1}, we can find a resulting CNF~\(S'\) equisatisfiable
  to~\(S\), and having no length~\(1\) clauses.  Since~\(S'\) has no
  length \(1\) clauses, it is a \mbox{\(2\)-CNF}.

  Having dealt with the case when~\(S\) has
  four~\mbox{\((a,b)\)-clauses} for some~\({a,b\in V}\), we can assume
  throughout the remainder of the proof that~\(S\) has three or
  fewer~\mbox{\((a,b)\)-clauses} for every~\({a,b\in V}\).  Also,
  having dealt with the case when~\(S\) has an~\mbox{\((a)\)-clause}
  for some~\(a\in V\), we can henceforth assume that~\(S\) is a
  reduced~\mbox{\(2\)-CNF}.

  Let~\(a,b\in V\) such that~\(S\) has exactly
  two~\mbox{\((a,b)\)-clauses}.  By exchanging~\(a\) with~\(\neg a\)
  \mbox{and/or}~\(b\) with~\(\neg b\) we can ensure that one of the
  two~\mbox{\((a,b)\)-clauses} is~\((a\vee b)\).  Thus~\(S\) has one
  of three forms, which we label~\(S_1, S_2\) and~\(S_3\)
  \begin{equation*}
    S_1 = R\wedge(a\vee b)\wedge(a\vee\neg b),
    \quad
    S_2 = R\wedge(a\vee b)\wedge(\neg a\vee b),
    \quad
    S_3 = R\wedge(a\vee b)\wedge(\neg a\vee\neg b);
  \end{equation*}
  where~\(R\) is a reduced~\mbox{\(2\)-CNF} not containing
  any~\mbox{\((a,b)\)-clauses}.  Lemma~\ref{lem:
    ab_clause_cases}~\eqref{item: ab_2_clause_1} implies
  that~\({S_1\sim R[a:=\top]}\).  We note that~\(R[a:=\top]\) is a
  reduced CNF with clauses of length either~\(2\) or~\(1\).  By
  repeated application of Lemma~\ref{lem: a_clause_cases}~\eqref{item:
    a_clause_2} and~\eqref{item: a_clause_1}, we can again find a
  reduced~\mbox{\(2\)-CNF}~\(R_1'\) such
  that~\(R_1'\sim R[a:=\top]\sim S_1\), and such that~\(R_1'\) (as
  opposed to~\(S_1\)) has no~\mbox{\((a,b)\)-clauses}.  Repeating this
  process for every~\(c,d\in V\) such that~\(R_1'\) has exactly
  two~\mbox{\((c,d)\)-clauses} gives us a
  resulting~\mbox{\(2\)-CNF}~\(S_1'\) equisatisfiable to~\(S_1\) such
  that~\(S_1'\) has either three, or one, or
  zero~\mbox{\((a,b)\)-clauses} for every~\(a,b\in V\).

  Lemma~\ref{lem: ab_clause_cases}~\eqref{item: ab_2_clause_2} implies
  that~\({S_2\sim R[b:=\top]}\).  Repeatedly applying
  Lemma~\ref{lem: a_clause_cases}~\eqref{item: a_clause_2}
  and~\eqref{item: a_clause_1}, we can find a
  reduced~\mbox{\(2\)-CNF}~\(R_2'\) such that~\(R_2'\sim S_2\) and
  such that~\(R_2'\) (as opposed to~\(S_2\)) has
  no~\mbox{\((a,b)\)-clauses}.  If \(R_2'\) has exactly
  two~\mbox{\((c,d)\)-clauses} for some~\(c,d\in V\), then repeating
  the above process gives us a resultant~\mbox{\(2\)-CNF}~\(S_2'\)
  equisatisfiable to~\(S_2\) such that~\(S_2'\) has either three, or
  one, or zero~\mbox{\((a,b)\)-clauses} for every~\(a,b\in V\).

  Lemma~\ref{lem: ab_clause_cases}~\eqref{item: ab_2_clause_3} implies
  that~\({S_3\sim R[b:=\neg a]}\).  By repeated application of
  Lemma~\ref{lem: a_clause_cases}~\eqref{item: a_clause_2}
  and~\eqref{item: a_clause_1}, we can again find an
  reduced~\mbox{\(2\)-CNF}~\(R_3'\) such that~\(R_3'\sim S_3\) and
  such that~\(R_3'\) (as opposed to~\(S_3\)) has
  no~\mbox{\((a,b)\)-clauses}.  If \(R_3'\) has exactly
  two~\mbox{\((c,d)\)-clauses} for some~\(c,d\in V\), then repeating
  the above process gives us a resultant~\mbox{\(2\)-CNF}~\(S_3'\)
  equisatisfiable to~\(S_3\) such that~\(S_3'\) has either four, or
  three, or one, or zero~\mbox{\((a,b)\)-clauses} for
  every~\(a,b\in V\).  If~\(S_3'\) has four~\mbox{\((a,b)\)-clauses}
  for any~\(a,b\in V\), then just as before, we can replace \(S_3'\)
  with \(\bot\).

  We can assume for the remainder of the proof that~\(S\) itself is a
  reduced \mbox{\(2\)-CNF} that has either three, or one, or
  zero~\mbox{\((a,b)\)-clauses} for every~\(a,b\in V\).
  Let~\({a,b\in V}\) be such that~\(S\) has exactly
  three~\mbox{\((a,b)\)-clauses}.  By exchanging~\(a\) with~\(\neg a\)
  \mbox{and/or}~\(b\) with~\(\neg b\) we can ensure that the
  three~\mbox{\((a,b)\)-clauses} are~\({(a\vee
    b)}\),~\({(a\vee\neg b)}\) and~\({(\neg a\vee b)}\).  Hence we can
  write~\({S = R\wedge(a\vee b)\wedge(a\vee\neg b)\wedge(\neg a\vee
    b)}\), where~\(R\) is a reduced~\mbox{\(2\)-CNF} not containing
  any~\mbox{\((a,b)\)-clauses}.  By Lemma~\ref{lem:
    ab_clause_cases}~\eqref{item: ab_3_clause}, we conclude
  that~\({S\sim R[\{a,b\}:=\top]}\).  The reduced
  CNF~\(R[\{a,b\}:=\top]\) does not have~\mbox{\((a,b)\)-clauses}.  By
  using procedures already described in this proof, we can find a
  reduced~\mbox{\(2\)-CNF}~\(R'\) equisatisfiable
  to~\(R[\{a,b\}:=\top]\) and having either three, or one, or
  zero~\mbox{\((c,d)\)-clauses} for every~\(c,d\in V\).  By repeating
  this procedure as many times as needed, we can find a resultant
  reduced~\mbox{\(2\)-CNF}~\(S'\) equisatisfiable to~\(S\) such
  that~\(S'\) does not have three~\mbox{\((a,b)\)-clauses} for
  any~\(a,b\in V\).  Therefore, the reduced~\mbox{\(2\)-CNF}~\(S'\)
  has one or fewer~\mbox{\((a,b)\)-clauses} for every~\(a,b\in V\).
  In other words, the reduced~\mbox{\(2\)-CNF}~\(S'\) is simple.
  \qed\end{pf}

Throughout the remainder of this paper, when we refer to graphs,
we will always mean simple graphs unless stated otherwise.

\section{\(\U\) is closed under graph homeomorphism}\label{sec:homeo}
In this section we look at three different graph operations --
subgraphing, subdivision and \mbox{edge-contraction} at edges not
contained in triangles -- and study if~\(\U\) is closed under these
operations.

A graph~\(G\) is a \emph{subgraph} of a graph~\(H\) if both the
vertex set and edge set of~\(G\) are subsets of the vertex and
edge sets of~\(H\).
We will show that if a subgraph of a graph is in~\(\U\),
then the graph itself must also be in~\(\U\).
\emph{Edge-contraction at an edge}~\((u,v)\) of a graph results
in a graph in which the vertices~\(u\) and~\(v\) are merged into a
single new vertex~\(w\) and the all edges incident on~\(u\) or~\(v\)
are now incident on~\(w\).
In order to avoid the creation of \mbox{multi-edges} we restrict the
\mbox{edge-contraction} operation to edges not contained in triangles.
This is a necessary restriction and without it one would lose the
correspondence between simple~\(2\)-CNFs and their associated simple
graphs.
We will show that if~\(G\) can be obtained from a graph in \(\U\) via
a series of \mbox{edge-contractions} at edges not contained in triangles,
then~\(G\) must be in~\(\U\).  

A \emph{subdivision of an edge}~\((u,v)\) in a graph yields a graph
containing one new vertex~\(w\), and with an edge set replacing~\((u,v)\)
by two new edges~\((u,w)\) and~\((w,v)\).
A \emph{subdivision of a graph}~\(G\) is a graph resulting from the
subdivision of edges in~\(G\).
Two graphs are \emph{homeomorphic} if they are subdivisions
of the same graph.
We will show that if two graphs are homeomorphic, then
either both are in~\(\U\) or neither is.

A graph~\(G\) is a \emph{topological minor} of a graph~\(H\) if a
subdivision of~\(G\) is a subgraph of~\(H\).
We will produce a complete list of graphs whose appearance as a
topological minor is an obstruction to satisfiability
of~\mbox{\(2\)-CNFs}.
It is possible to embed simple graphs into~\(\mathbb{R}^3\) and
to allow them to inherit the subspace topology of~\(\mathbb{R}^3\).
Once embedded, homeomorphic graphs are also homeomorphic in the
topological sense and topological graph minors are simply
topological subspaces.
In that sense, the ability to support unsatisfiable sentences is a
topological property of graphs.

\vspace*{3mm}
\begin{thm} \label{thm: U_closure} If simple graphs \(G\) and \(H\)  are
  homeomorphic, then~\({G\in\U}\) if and only if~\({H\in\U}\).
\end{thm}
\begin{pf}
  Since~\(G\) and~\(H\) are homeomorphic graphs, there exists a
  graph~\(K\) such that both~\(G\) and~\(H\) are subdivisions of~\(K\).
  It is enough to prove that~\({K\in\U}\) if and only if~\({G\in\U}\).
  In fact, it is enough to prove that~\({K\in\U}\) if and only
  if~\({G'\in\U}\), where~\(G'\) is a graph obtained via a single
  subdivision at an arbitrary edge~\((u,v)\) in~\(K\).
  We denote by~\(w\) the new vertex in~\(G'\) created by the subdivision.
  
  Suppose that~\({K\in\U}\).
  Then, there exists an unsatisfiable simple~\(2\)-CNF~\(S_K\) supported
  on~\(K\).
  We can write without loss of generality,~\(S_K = S\wedge (u\vee v)\),
  where~\(S\) is a simple~\mbox{\(2\)-CNF} not
  containing~\mbox{\((u,v)\)-clauses}.
  Note that if~\(S_K\) is not in this form, that is, if
  the~\mbox{\((u,v)\)-clause} is not positive in either~\(u\) or~\(v\),
  then we simply interchange~\(u\) with~\(\neg u\) \mbox{and/or}~\(v\)
  with~\(\neg v\) in order to bring it to the desired form.
  This interchange of a literal with its negation results in a different,
  but equisatisfiable simple~\(2\)-CNF.
  By abuse of notation, we will continue to refer to this new~\(2\)-CNF
  as~\(S_K\).
  We then define a simple~\mbox{\(2\)-CNF}~\(S_{G'}\) supported on \(G'\)
  as \(S_{G'} = S \wedge \paren*{u\vee w} \wedge\paren*{\neg w\vee v}\).
  We prove that~\(S_{G'}\) is unsatisfiable by showing that from any
  truth assignment that satisfies~\(S_{G'}\), we can obtain a truth
  assignment satisfying~\(S_K\), leading to a contradiction.
  
  In any truth assignment for~\(S_{G'}\), either~\(w\) is set to~\(\top\)
  or it is set to~\(\bot\).
  If~\([w:=\top]\), then it implies
  \begin{equation*}
    S_{G'}[w:=\top] \;=\; S \wedge v \;\sim\; S[v:=\top] \;=\; S_K[v:=\top],
  \end{equation*}
  is satisfiable, that is, the~\mbox{\(2\)-CNF}~\(S_K\) is satisfiable by
  some truth assignment which sets~\({[v:=\top]}\).
  Similarly, if \([w:=\bot]\), then it implies
  \begin{equation*}
  S_{G'}[w:=\bot] \;=\; S \wedge u \;\sim\; S[u:=\top] \;=\; S_K[u:=\top],
  \end{equation*}
  is satisfiable, that is, the~\mbox{\(2\)-CNF}~\(S_K\) is satisfiable by some truth assignment that
  sets~\({[u:=\top]}\).
  Either scenario leads to a contradiction since~\(S_K\) is
  assumed to be unsatisfiable.
  Thus, we conclude that~\(S_{G'}\) is also unsatisfiable and
  that~\({G'\in\U}\).
  
  Conversely, suppose that~\({G'\in\U}\).
  Then, there exists an unsatisfiable simple~\mbox{\(2\)-CNF}~\(S_{G'}\),
  supported on~\(G'\).
  By exchanging~\(u\) with~\(\neg u\), \mbox{and/or}~\(v\)
  with~\(\neg v\) \mbox{and/or}~\(w\) with~\(\neg w\), we can assume
  without loss of generality that either~\(S_{G'}\) has the form
  \begin{equation*}
    S_{G'} = S \wedge \paren*{u\vee w} \wedge\paren*{w\vee v}
    \qquad\text{or}\qquad
    S_{G'} = S \wedge \paren*{u\vee w} \wedge\paren*{\neg w\vee v}.
  \end{equation*}
  In the first case, since~\(S_{G'}\) is unsatisfiable,
  by setting~\({[w:=\top]}\) we obtain an unsatisfiable~\mbox{\(2\)-CNF}.
  Thus \({S_{G'}[w:=\top] = S}\) is unsatisfiable.
  Therefore, in this case, we have found~\({S\wedge (u\vee v)}\), an
  unsatisfiable~\mbox{\(2\)-CNF} supported on~\(K\).
  In the other case, since~\(S_{G'}\) is unsatisfiable,
  the~\mbox{\(2\)-CNF} obtained by setting~\([w:=v]\) must also be
  unsatisfiable.
  Thus,~\({S_{G'}[w:=v] = S\wedge(u\vee v)}\) is an
  unsatisfiable~\mbox{\(2\)-CNF} supported on~\(K\).
  We conclude that~\({K\in\U}\).\qed\end{pf}

In view of the theorem we just proved, it is natural to attempt a
classification of all members of~\(\U\) up to homeomorphism.
However, we state here (without proof) that there are infinitely many
graphs in~\(\U\) even up to homeomorphism.
Given below is a (\mbox{non-exhaustive}) infinite list of graphs
in~\(\U\), none of which are homeomorphic to each other.
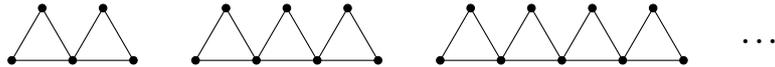
\[\Twohills  \qquad \Threehills  \qquad \Fourhills \qquad \boldsymbol{\ldots}\]

\begin{thm}\label{thm: subgraph_U}
  Let~\(G\) and~\(H\) be simple graphs such that~\(G\) is a subgraph
  of~\(H\).
  If~\({G\in\U}\), then~\({H\in\U}\).
\end{thm}
\begin{pf}
  If \(G\in\U\) then there exists an unsatisfiable simple
  \mbox{\(2\)-CNF} \(S_G\) supported on~\(G\).
  We denote by \(S\) a simple \mbox{\(2\)-CNF} supported on the
  graph with vertex set~\(V(H)\) and edge set~\(E(H)-E(G)\).
  We then define a simple~\mbox{\(2\)-CNF}~\(S_H\)
  as~\({S_H=S\wedge S_G}\).
  The \mbox{\(2\)-CNF}~\(S_H\) is supported on~\(H\) and is
  unsatisfiable since any truth assignment satisfying~\(S_H\)
  would also satisfy~\(S_G\).
  Hence~\(H\in\U\).\qed\end{pf}

We note that the converse of this theorem is not true, that is,
there exist graphs~\(G\notin\U\) and~\(H\in\U\) such that~\(G\)
is a subgraph of~\(H\).
For example, the triangle graph is not in~\(\U\) (as shown in
Lemma~\ref{lem: C3_is_sat}) but~\(\;\;\WhiteButterfly\in\U\) (as shown
in the proof of Lemma~\ref{lem: minimal_unsat_graphs}).

\begin{cor}\label{cor: topological_minor_in_U}
  Let simple graph~\(G\) be a topological minor of a simple graph~\(H\).
  If~\(G\in\U\) then \(H\in\U\).
\end{cor}
\begin{pf}
  By definition of topological minors, some subdivision~\(G'\)
  of~\(G\) is isomorphic to a subgraph of~\(H\).
  By Theorem~\ref{thm: U_closure}, if~\(G\in\U\) then~\(G'\in\U\).
  By Theorem~\ref{thm: subgraph_U}, if~\(G'\in\U\) then~\(H\in\U\).
  \qed
\end{pf}

We note here that the converse of
Corollary~\ref{cor: topological_minor_in_U} is not true for the
same reason that the converse of Theorem~\ref{thm: subgraph_U}
is not.
Corollary~\ref{cor: topological_minor_in_U} motivates the
following definition --- a graph~\(G\) is a \emph{minimal unsatisfiability
  graph} if both of the following conditions hold
\begin{enumerate}[a)]
  \item \(G\in\U\),
  \item for every proper topological minor~\(G'\) of~\(G\), we have
  that~\({G'\notin\U}\).
\end{enumerate}
Furthermore, a set~\(M\) of minimal unsatisfiability graphs is
\emph{complete} if every graph in~\(\U\) has a subgraph that is 
homeomorphic to some element of~\(M\), that is, if every graph in~\(\U\)
has an element of~\(M\) as a topological minor.

Since we know that minimal unsatisfiability graphs exist, we can form a
set~\(M\) by constructing the union of all sets of minimal unsatisfiability
graphs. This set is complete because each graph in~\(\U\) will have
some element of~\(M\) as a topological minor -- if not we simply add
this to~\(M\) a minimal unsatisfiability graph that is a topological
minor of this graph. Further, such a complete set must be unique
because if not, then there is some minimal unsatisfiability graph~\(G\)
in complete set~\(M'\) that is not in~\(M\). This would result in
\(G\) having a proper topological minor in~\(\U\), violating the
minimality of \(G\).
After proving the following result about the relation between
\mbox{edge-contraction} and graphs in~\(\U\), the remainder of this paper
is dedicated to finding this unique complete set of minimal unsatisfiability graphs.

\begin{thm}\label{thm: edge_contraction_U}
  Let \(G\) and \(H\) be simple graphs such that~\(G\) can be obtained
  via a series of \mbox{edge-contractions} at edges of~\(H\) not
  contained in triangles.
  If~\({H\in\U}\), then~\({G\in\U}\).
\end{thm}
\begin{pf}
  If is enough to prove the theorem for the case when~\(G\) can be obtained
  from \(H\) via a single \mbox{edge-contraction}, say at the edge~\((u,v)\)
  not contained in any triangles in~\(H\).
  We label the new vertex in~\(G\) formed by the merger of~\(u\) and~\(v\)
  by~\(w\).
  
  Suppose that \({H\in\U}\).
  Then, there exists an unsatisfiable simple~\mbox{\(2\)-CNF}~\(S_H\),
  supported on~\(H\).
  Without loss of generality, we can write
  \begin{equation*}
    S_H = S_1 \wedge (u\vee v)
  \end{equation*}
  where~\(S_1\) is a simple~\mbox{\(2\)-CNF} not containing
  any~\mbox{\((u,v)\)-clauses}. 
  If~\(S_H\) is not in this form, we can exchange~\(u\) with~\(\neg u\)
  \mbox{and/or}~\(v\) with~\(\neg v\) to make the~\mbox{\((u,v)\)-clause}
  in~\(S_H\) positive in both variables.
  
  One can further factor~\(S_1\) such that
  \begin{equation*}
  S_H = S_2 \wedge \paren*{\bigwedge_{a\in A} a\vee u}
            \wedge \paren*{\bigwedge_{b\in B} b\vee\neg u}
            \wedge \paren*{u\vee v},
  \end{equation*}
  where~\(S_2\) is a simple~\(2\)-CNF not containing any clauses incident
  on the vertex~\(u\), and the sets \(A\) and \(B\) are disjoint subsets
  of~\(\paren*{V\cup\neg V}-\{u, \neg u\}\).
  The sets~\(A\) and~\(B\) are disjoint because~\(S_H\) is a simple~\(2\)-CNF.
  
  Next, we factor~\(S_2\) such that
  \begin{equation*}
    S_H = S_3 \wedge \paren*{\bigwedge_{a\in A} a\vee u}
              \wedge \paren*{\bigwedge_{b\in B} b\vee\neg u}
              \wedge \paren*{\bigwedge_{c\in C} c\vee v}
              \wedge \paren*{\bigwedge_{d\in D} d\vee\neg v}
              \wedge \paren*{u\vee v},
  \end{equation*}
  where~\(S_3\) is a simple~\(2\)-CNF not containing any clauses incident
  on vertices~\(u\) or~\(v\), and the sets \(A,B,C,D\) are \mbox{pairwise-disjoint}
  subsets of~\(\paren*{V\cup\neg V}-\{u, \neg u, v, \neg v\}\).
  The sets~\(C\) and~\(D\) are disjoint because~\(S_H\) is a simple~\(2\)-CNF.
  The sets \(A\) and~\(C\) are disjoint since the edge~\((u,v)\) is not contained
  in a triangle.
  Disjointness of other pairs of sets follows similarly.
  
  We then choose
  \[S_G = S \wedge \paren*{\bigwedge_{x\in A\cup D} x\vee w}
            \wedge \paren*{\bigwedge_{y\in B\cup C} y\vee\neg w},\]
  and note that~\(S_G\) is a simple \mbox{\(2\)-CNF} supported on~\(G\).
  The pairs of sets~\((A,D)\) and~\((B,C)\) being disjoint implies that no clauses are
  lost in the process of \mbox{edge-contraction}. 
  The pairs~\((A,B)\),~\((A,C)\),~\((B,D)\) and~\((C,D)\) being disjoint implies
  that~\(S_G\) is simple.
  
  We prove that~\(S_G\) is unsatisfiable by showing that from any
  truth assignment that satisfies \(S_G\), we can obtain a truth
  assignment satisfying~\(S_H\), leading to a contradiction.
  
  Given a truth assignment for \(S_G\), we can extend it to a truth
  assignment for~\(S_H\) by setting~\({[u:=w]}\) and~\({[v:=\neg w]}\).
  We have
  \begin{equation*}
    S_H[u:=w][v:=\neg w] \;=\;
       S\wedge\paren*{\bigwedge_{x\in A\cup D} x\vee w}
        \wedge \paren*{\bigwedge_{y\in B\cup C} y\vee\neg w} \;=\; S_G,
  \end{equation*}
  is satisfiable.
  This contradicts the unsatisfiability of \(S_H\).
  We conclude that \(S_G\) is unsatisfiable, and hence that~\({G\in\U}\).
  \qed\end{pf}

The converse of this theorem is not true, that is, there
does exist a graph~\(G\in\U\) that can be obtained via \mbox{edge-contractions}
of a graph~\(H\notin\U\) at edges not contained in triangles of~\(H\).
For example, consider the graphs~\({G =\;\WhiteButterfly}\) and
\({H=\;\SquareButterfly}\).
The graph~\(G\) is in~\(\U\) (as shown in
Lemma~\ref{lem: minimal_unsat_graphs}) and can be obtained from an
\mbox{edge-contraction} from~\(H\). However, the graph~\(H\) is not in~\(\U\)
(implied by Theorem~\ref{thm: edge_contraction_U} and Lemma~\ref{lem: K4-e_is_sat}
when combined with the fact that~\(H\) can be reduced
to~\({K_4-e}\) via other \mbox{edge-contractions}).

\section{Graph families that are not in \(\U\)}

\begin{lem}\label{lem: C3_is_sat}
  Let~\(C_3\) denote the triangle graph.
  The graph~\(C_3\) is not in~\(\U\).
\end{lem}
\begin{pf}
  We enumerate the vertex set of~\(C_3\) as~\(\{a,b,c\}\).
  Without loss of generality, every simple~\mbox{\(2\)-CNF}
  supported on~\(C_3\) can be written either in the form
  \begin{equation*}
    S = (a\vee b) \wedge (a \vee c) \wedge (b,c)\text{-clause}
    \qquad\text{or}\qquad
    S = (a\vee b) \wedge (\neg a \vee c) \wedge (b,c)\text{-clause}.
  \end{equation*}
  If~\(S\) is not originally in this form, we modify it by
  interchanging~\(a\) with~\(\neg a\) \mbox{and/or}~\(b\)
  with~\(\neg b\) \mbox{and/or}~\(c\) with~\(\neg c\) till it is.
  
  In the first case, setting~\([a:=\top]\) gives yields a
  single~\mbox{\((b,c)\)-clause}.
  This clause can be satisfied by making the appropriate
  assignment for either~\(b\) or~\(c\).
  In the second case, by setting \([\{a, c\}:=\top]\) we get
  \begin{equation*}
    S[\{a, c\}:=\top] \;=\; (b,c)\text{-clause}\,[c:=\top]
    \;\sim\;(b,c)\text{-clause}\wedge c.
  \end{equation*}
  The~\mbox{\((c)\)-clause} can be satisfied by
  setting~\([c:=\top]\), while~\(b\) can be set to an appropriate
  assignment in order to satisfy the~\mbox{\((b,c)\)-clause}.
  Thus, the~\mbox{\(2\)-CNF}~\(S\) is satisfiable.
  \qed\end{pf}

\begin{cor}
  Let \(C_n\) denote the cycle graph on~\(n\) vertices.
  The graph~\(C_n\) is not in~\(\U\) for any~\(n\geq 3\). 
\end{cor}
\begin{pf}
  The graph~\(C_n\) is homeomorphic to~\(C_3\) for
  every~\(n\geq 3\).
  Theorem~\ref{thm: U_closure} and Lemma~\ref{lem: C3_is_sat}
  therefore imply that \(C_n\notin\U\).\qed
\end{pf}

\begin{lem}\label{lem: K4-e_is_sat}
  Let~\(K_4\) denote the complete graph on four vertices.
  Let~\({K_4-e}\) denote the graph obtained by deleting a single
  edge~\(e\) from~\(K_4\).
  The graph~\(K_4-e\) is not in~\(\U\).
\end{lem}
\begin{pf}
  We enumerate the vertex set of~\(K_4\) as~\(\{a,b,c,d\}\).
  Let~\(e\) denote the edge~\({(c,d)\in E(K_4)}\) such
  that~\({(c,d)\notin E(K_4-e)}\).
  Every~\mbox{\(2\)-CNF}~\(S\) supported on~\({K_4-e}\) can be
  written either in the form
  \begin{align*}
    S &= (a\vee b)\wedge(a\vee c)\wedge(a\vee d)\wedge
         (b,c)\text{-clause}\wedge(b,d)\text{-clause},
    \quad\text{or}\\
    S &= (a\vee b)\wedge(a\vee c)\wedge(\neg a\vee d)
         \wedge (b,c)\text{-clause}\wedge(b,d)\text{-clause}.
  \end{align*}
  If~\(S\) is not already in the desired form, then we can
  interchange each variable with its negation till it is.
  
  In the first case, by setting~\([a:=\top]\) we
  obtain
  \begin{equation*}
    {S[a:=\top] = (b,c)\text{-clause}\wedge(b,d)
    \text{-clause}}.
  \end{equation*}
  This can be satisfied by making appropriate assignments
  for~\(c\) and~\(d\) so that they satisfy each of the clauses.
  In the second case, we can set~\([a:=\top]\) to
  obtain
  \begin{equation*}
    {S[a:=\top] = d\wedge(b,c)\text{-clause}\wedge(b,d)
      \text{-clause}}.
  \end{equation*}
  This resulting CNF can be satisfied by setting~\([d:=\top]\),
  by choosing an assignment for~\(b\) that would satisfy
  the~\mbox{\((b,d)\)-clause} and by then choosing an assignment
  for~\(c\) that would satisfy the~\mbox{\((b,c)\)-clause}.
  We conclude that~\(S\) is always satisfiable, and therefore,
  the graph~\(K_4-e\) is not in~\(\U\).
  \qed\end{pf}

\begin{lem}\label{lem: trees_are_sat}
  Tree graphs are not in \(\U\).
\end{lem}
\begin{pf}
  Every tree graph can be reduced via \mbox{edge-contractions} to
  a \mbox{single-vertex} graph.
  A \mbox{single-vertex} graph is clearly not in \(\U\).
  The result follows from Theorem~\ref{thm: edge_contraction_U}.\qed
\end{pf}

\section{Structure of graphs with two or three independent cycles} 
\label{sec:structure}

In this section we prove four lemmas about the structure of graphs that
have either two or three independent cycles.
These structural results are needed for proving the results in
Section~\ref{sec:complete}, including the main result of this
paper~(Theorem~\ref{thm:main}).

\begin{lem}\label{lem: two_C3_combinations}
  Every connected simple graph having exactly two copies of~\(C_3\)
  as subgraphs has one of the following three graphs as a topological
  minor.
  
  \centering
  {\setlength{\tabcolsep}{2em}
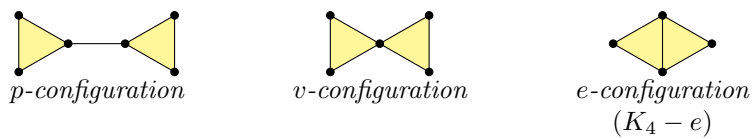
    \begin{tabular}{ccc}
      \Pconfig            & \Vconfig            & \Econfig\\
      \(p\)-configuration & \(v\)-configuration & \(e\)-configuration\\
                          &                     & \((K_4-e)\)
  \end{tabular}}
\end{lem}

\begin{rmk}
  For convenience, we label the three graphs as the~\mbox{\(p\)-configuration}
  (\(p\)~stands for path), the
  \mbox{\(v\)-configuration} (\(v\) stands for vertex) and
  \mbox{\(e\)-configuration} (\(e\)~stands for edge) respectively.
  The two copies of~\(C_3\) have been filled in for easy visual
  identification.
\end{rmk}

\begin{pf}
  We construct this list of topological minors from the bottom
  up.
  Two copies of~\(C_3\) can be put together to create a connected simple
  graph in exactly three ways
  \begin{enumerate}[a)]
    \item either the two copies of~\(C_3\) share zero vertices,
    \item or they share exactly one vertex,
    \item or they share exactly two vertices (that is, they share
    an edge).
  \end{enumerate}
  In the first case, since the graph is still supposed to be
  connected, we claim that the two copies of~\(C_3\) must be
  connected by one or more paths.
  Such a graph will always have the two copies connected by a
  single path as a subgraph.
  Hence, the graph will also always have
  the~\mbox{\(p\)-configuration} as a topological minor.
  
  In the second case, the graph will always have the
  \mbox{\(v\)-configuration} as a subgraph and therefore, also as a
  topological minor.
  Similarly, in the third case, the graph will always have
  the~\mbox{\(e\)-configuration} as a subgraph and therefore, also as
  as topological minor.
  \qed\end{pf}

\begin{lem}\label{lem: two_cycle_combinations}
  Every connected simple graph having exactly two or more independent
  cycles has one of the three graphs in Lemma~\ref{lem: two_C3_combinations}
  as a topological minor.
\end{lem}
\begin{pf}
  Let~\(G\) be a connected simple graph with two or more independent cycles.
  If any two cycles share one or more edges, then~\(K_4-e\) is a
  topological minor of~\(G\).
  If no pair of cycles share any edges, but at least one pair shares a vertex,
  then the~\mbox{\(v\)-configuration} is a topological minor of~\(G\).
  If every pair of cycles share neither any edges nor any vertices, then using
  the connectedness of~\(G\), we infer that there must be a path connecting
  vertices in every pair of cycles.
  Thus, in this case, the~\mbox{\(p\)-configuration} is a topological minor
  of~\(G\).\qed
\end{pf}

\begin{lem}\label{lem: three_triangle_combinations}
  Every connected simple graph having three or more copies of~\(C_3\)
  as subgraphs has one of the following~\(15\) graphs as a topological
  minor.
  \newcommand{\vspacing}{\vspace*{3ex}}
  
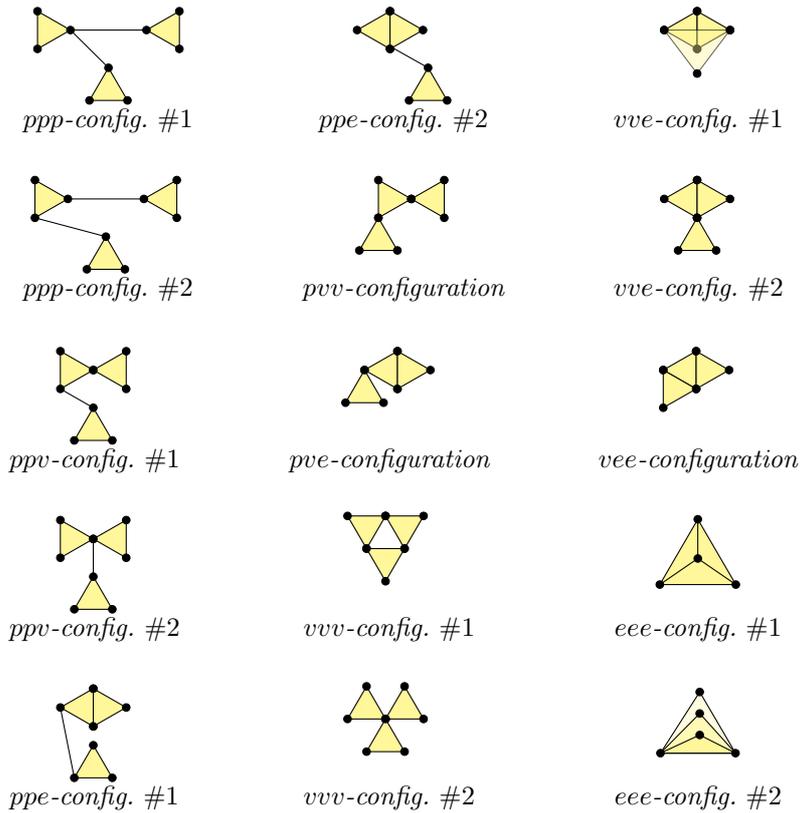
  \begin{center}
  \setlength{\tabcolsep}{2em}
  \begin{tabular}{ccc}
    \PPPoneconfig            & \PPEtwoconfig            & \VVEoneconfig\\
    \(ppp\)-config.\ \(\#1\) & \(ppe\)-config.\ \(\#2\) & \(vve\)-config.\ \(\#1\)\vspacing\\

    \PPPtwoconfig            & \PVVconfig               & \VVEtwoconfig\\
    \(ppp\)-config.\ \(\#2\) & \(pvv\)-configuration    & \(vve\)-config.\ \(\#2\)\vspacing\\
  \end{tabular}
  
  \begin{tabular}{ccc}
    \PPVoneconfig            & \PVEconfig               & \VEEconfig\\
    \(ppv\)-config.\ \(\#1\) & \(pve\)-configuration    & \(vee\)-configuration\vspacing\\
    
    \PPVtwoconfig            & \VVVoneconfig            & \EEEoneconfig\\
    \(ppv\)-config.\ \(\#2\) & \(vvv\)-config.\ \(\#1\) & \(eee\)-config.\ \(\#1\)\vspacing\\
    
    \PPEoneconfig            & \VVVtwoconfig            & \EEEtwoconfig\\
    \(ppe\)-config.\ \(\#1\) & \(vvv\)-config.\ \(\#2\) & \(eee\)-config.\ \(\#2\)
  \end{tabular}
  \end{center}
\end{lem}
\begin{rmk}
  Labels for each of the~\(15\) graphs
  (like~\mbox{\(ppp\)-config.}~\(\#1\)) are explained as
  part of the proof.
  The three copies of~\(C_3\) that we use for this labeling
  have been filled in for the sake of easy visual identification.
\end{rmk}
\begin{pf}
  We arrive at this list of topological minors by constructing them
  from the bottom up.
  From the proof of Lemma~\ref{lem: two_C3_combinations}, we
  know that every pair of~\(C_3\) can be joined in exactly
  one of three ways and that these three ways are inequivalent
  when considering the topological minor relation, that is, no graph
  in Lemma~\ref{lem: two_C3_combinations} is a topological minor of
  another graph in Lemma~\ref{lem: two_C3_combinations}. 
  Graphs with three or more copies of~\(C_3\) will have at
  least~\({3\choose 2}=3\) distinct pairs of~\(C_3\).
  
  To enumerate all possible joinings of these three pairs of~\(C_3\), we form
  all \mbox{three-letter} words formed using the
  letters~\(\{p, v, e\}\), with repetition allowed, but order being
  irrelevant.
  For example, the word~\(ppe\) would correspond to the family of
  graphs where the two pairs of~\(C_3\) are joined via
  the~\mbox{\(p\)-configuration}, while the third pair
  of~\(C_3\) is joined via the~\mbox{\(e\)-configuration}.
  We also note that some words correspond to multiple
  \mbox{non-isomorphic} configurations.
  We analyze each case separately below
  \begin{enumerate}[a)]
    \item \(ppp\)-configuration. After joining the first pair
    of~\(C_3\) in a~\mbox{\(p\)-configuration},
    the third copy of~\(C_3\) can be added in two different
    \mbox{non-isomorphic} ways (such that the first and third,
    as well as the second and third are joined via
    the~\mbox{\(p\)-configuration}).
    Thus, we obtain only two configurations in this
    case~---~\(ppp\)-config.~\(\#1\) and~\(ppp\)-config.~\(\#2\).
    
  \item \(ppv\)-configuration.
    After joining the first pair of~\(C_3\) in
    a~\mbox{\(v\)-configuration}, the third copy of~\(C_3\) can be
    added in two \mbox{non-isomorphic} ways.
    Thus, we obtain only two configurations in this
    case, namely~\(ppv\)-config.~\(\#1\) and~\(ppv\)-config.~\(\#2\).
    
  \item \(ppe\)-configuration.
    After joining the first pair of~\(C_3\) in
    an~\mbox{\(e\)-configuration}, the third copy of~\(C_3\) can be
    added in two \mbox{non-isomorphic} ways.
    Thus, we obtain only two configurations in this
    case, namely \(ppe\)-config.~\(\#1\) and~\(ppe\)-config.~\(\#2\).
  
  \item \(pvv\)-configuration.
    After joining the first pair of~\(C_3\) in
    a~\mbox{\(v\)-configuration}, the third copy of~\(C_3\) can be
    added in only one way.
    Thus, we obtain only one configurations in this
    case, namely the~\(pvv\)-configuration.
  
  \item \(pve\)-configuration.
    After joining the first pair of~\(C_3\) in
    an~\mbox{\(e\)-configuration}, the third copy of~\(C_3\) can
    added in only one way.
    Thus, we obtain only one configurations in this
    case, namely the~\(pve\)-configuration.
    
  \item \(pee\)-configuration.
    After joining the first pair of~\(C_3\) in
    an~\mbox{\(e\)-configuration}, the third copy of~\(C_3\) cannot
    be added in a way that it is shares an edge with the first and
    is connected by a path to the second.
    Thus, we do not obtain configurations in this case.

  \item \(vvv\)-configuration.
    After joining the first pair of~\(C_3\) in
    an~\mbox{\(v\)-configuration}, the third copy of~\(C_3\) can be
    added in two \mbox{non-isomorphic} ways.
    Thus, we obtain only two configurations in this
    case, namely \(vvv\)-config.~\(\#1\) and~\(vvv\)-config.~\(\#2\).
  
  \item \(vve\)-configuration.
    After joining the first pair of~\(C_3\) in
    an~\mbox{\(e\)-configuration}, the third copy of~\(C_3\) can be
    added in two \mbox{non-isomorphic} ways.
    Thus, we obtain only two configurations in this
    case, namely \(vve\)-config.~\(\#1\) and~\(vve\)-config.~\(\#2\).

  \item \(vee\)-configuration.
    After joining the first pair of~\(C_3\) in
    an~\mbox{\(e\)-configuration}, the third copy of~\(C_3\) can
    added in only one way.
    Thus, we obtain only one configurations in this
    case, namely the~\(vee\)-configuration.

  \item \(eee\)-configuration.
    After joining the first pair of~\(C_3\) in
    an~\mbox{\(e\)-configuration}, the third copy of~\(C_3\) can be
    added in two \mbox{non-isomorphic} ways.
    Thus, we obtain only two configurations in this
    case, namely \(eee\)-config.~\(\#1\) and~\(eee\)-config.~\(\#2\).\qed
  \end{enumerate}
  \end{pf}

\begin{lem}\label{lem: three_cycle_combinations}
  Every connected simple graph having three or more independent cycles
  has one of the following four graphs as a topological minor.
  
  \centering
  {\setlength{\tabcolsep}{1em}
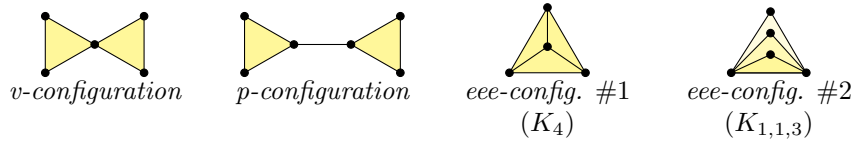
    \begin{tabular}{cccc}
      \Vconfig                 & \Pconfig                 & \EEEoneconfig                & \EEEtwoconfig\\
      \textit{v-configuration} & \textit{p-configuration} & \textit{eee-config. }\(\#1\) & \textit{eee-config. }\(\#2\)\\
                               &                          & \((K_4)\)                    & \((K_{1,1,3})\)
    \end{tabular}}

\end{lem}
\begin{pf}
  It is enough to show that every graph listed in
  Lemma~\ref{lem: three_triangle_combinations} has one of the
  four graphs listed above as a topological minor.
  The~\mbox{\(p\)-config.} is a subgraph of
  both the~\mbox{\(ppp\)-configs.}, both
  the~\mbox{\(ppv\)-configs.}, both~\mbox{\(ppe\)-configs.},
  the~\mbox{\(pvv\)-config.}, and the~\mbox{\(pve\)-config}.
  Of the remaining configurations,
  both the~\mbox{\(vvv\)-configs.}, both the~\mbox{\(vve\)-configs.},
  and the~\mbox{\(vee\)-config.} contain the~\mbox{\(v\)-config.}
  as a subgraph.
  
  The first~\mbox{\(eee\)-configuration} is isomorphic
  to~\(K_4\) while the second is isomorphic to \(K_{1,1,3}\).
  Hence, all~\(15\) configurations have at least one of the
  four graphs listed above as a topological minor.\qed\end{pf}

\section{The complete set of minimal unsatisfiability graphs}
\label{sec:complete}

We note that a simple graph~\(G\) is in~\(\U\) if and only if some
connected component of~\(G\) is in~\(\U\).  We proceed to make some
more observations about the \mbox{graph-family}~\(\U\).

\begin{lem}\label{lem: minimal_unsat_graphs}
  The four graphs in Lemma~\ref{lem: three_cycle_combinations} are
  minimal unsatisfiability graphs.
\end{lem}
\begin{pf}
  Consider the following unsatisfiable~\mbox{\(2\)-CNFs}.
  \begin{align*}
    S_1 &= (a\vee b)
            \wedge(\neg a\vee c)
            \wedge(\neg b\vee c)
            \wedge(\neg c\vee d)
            \wedge(\neg c\vee e)
            \wedge(\neg d\vee\neg e), \\
    S_2 &= (a\vee b)
           \wedge(\neg a\vee c)
           \wedge(\neg b\vee c)
           \wedge(\neg c\vee d)
           \wedge(\neg d\vee e)
           \wedge(\neg d\vee f)
           \wedge(\neg e\vee \neg f),\\
    S_3 &= (a\vee b)
           \wedge(a\vee c)
           \wedge(\neg a\vee d)
           \wedge(\neg b\vee\neg c)
           \wedge(b\vee\neg d)
           \wedge(c\vee\neg d)\;\mbox{ and} \\
    S_4 &= (a\vee b)
           \wedge(\neg a\vee d)
           \wedge(b\vee c)
           \wedge(\neg b\vee d)
           \wedge(\neg b\vee e)
           \wedge(\neg c\vee\neg d)
           \wedge(\neg d\vee\neg e).
  \end{align*}
  Their associated graphs are the four graphs listed above in order.
  Since these \mbox{\(2\)-CNFs} are unsatisfiable, the graphs are
  in~\(\U\).
  To prove that they are minimal unsatisfiability graphs, we need
  to show that none of the proper topological minors of these graphs
  are in~\(\U\). 
  Since the four listed graphs are not subdivisions of
  other graphs, it is enough to show that every proper
  subgraph of these graphs is not in~\(\U\).
    
  Every proper subgraph of the \mbox{\(v\)-configuration}
  is either a subgraph of \ButterflySubgraphOne{} or
  \ButterflySubgraphTwo{}.
  Both of these graphs can be reduced via
  \mbox{edge-contractions} to~\(C_3\).
  In Lemma~\ref{lem: C3_is_sat} we proved~\mbox{\(C_3\notin\U\)}
  and hence by Theorem~\ref{thm: edge_contraction_U} we
  can conclude that neither \ButterflySubgraphOne{} nor
  \ButterflySubgraphTwo{} is in \(\U\).
  Therefore, by Theorem~\ref{thm: subgraph_U}, none of the
  subgraphs of~\(v\)-configuration are in~\(\U\).
  
  Every proper subgraph of the \mbox{\(p\)-configuration}
  is either a subgraph of \BowtieSubgraphOne{},
  \BowtieSubgraphTwo{} or \BowtieSubgraphThree{}.
  The first two of these graphs can be reduced via
  edge-contractions to proper subgraphs of
  the~\mbox{\(v\)-configuration}.
  Since we proved that none of the proper subgraphs of
  the~\mbox{\(v\)-configuration} are in~\(\U\),
  by Theorem~\ref{thm: edge_contraction_U} we
  can conclude that \BowtieSubgraphOne{} nor
  \BowtieSubgraphTwo{} is in~\(\U\).
  Lastly, the graph~\BowtieSubgraphThree{} is in~\(\U\) if and
  only if at least one of its connected components is
  in \(\U\).
  In light of Lemma~\ref{lem: C3_is_sat}, we conclude that
  \BowtieSubgraphThree{} is not in~\(\U\).
  Theorem~\ref{thm: subgraph_U} implies that
  none of the proper subgraphs of
  the~\(p\)-configuration are in~\(\U\).
  
  In Lemma~\ref{lem: K4-e_is_sat}, we
  proved that~\(K_4-e\notin\U\).
  Since every proper subgraph of~\(K_4\) is also a subgraph
  of~\(K_4-e\), using Theorem~\ref{thm: subgraph_U} we
  conclude that none of the proper subgraphs of~\(K_4\) are
  in~\(\U\).
  
  Every proper subgraph of~\(K_{1,1,3}\) is either a
  subgraph \mbox{of~\;\BookSubgraphOne} or \(K_{2,3}\)~(~\BookSubgraphTwo),
  both of which can be reduced via
  \mbox{edge-contractions} to~\(K_4-e\).
  In Lemma~\ref{lem: K4-e_is_sat} we proved~\({K_4-e\notin\U}\)
  and hence by Theorem~\ref{thm: edge_contraction_U} we
  conclude that none of the proper subgraphs of~\(K_{1,1,3}\)
  are in~\(\U\).
  \qed\end{pf}

\begin{cor}
  Every connected simple graph having three or more independent
  cycles is in~\(\U\).
\end{cor}
\begin{pf}
  From Lemma~\ref{lem: three_cycle_combinations} we know that
  every connected simple graph having three or more independent
  cycles has one of the four graphs listed in that lemma as a
  topological minor.
  Since we proved in
  Lemma~\ref{lem: minimal_unsat_graphs} that all four
  of these graphs are in \(\U\), the result follows from
  Corollary~\ref{cor: topological_minor_in_U}.\qed\end{pf}

\begin{thm}\label{thm:main}
  The set~\({\raisebox{2.5pt}{\Big\{}\;\Butterfly,\;\;\Bowtiesmall,\;\;
    \Kfour,\;\;\Book\raisebox{2.5pt}{\Big\}}}\) is the complete set
  of minimal unsatisfiability graphs.
\end{thm}
\begin{pf}
  We denote the set by \(M\).
  Since we proved in Lemma~\ref{lem: minimal_unsat_graphs} that
  the elements of~\(M\)
  are minimal unsatisfiability graphs, it suffices to show that
  a simple graph~\(G\) is in~\(\U\) if and only if~\(G\) has some
  element of~\(M\) as a topological minor.
  Corollary~\ref{cor: topological_minor_in_U} implies the ``if''
  part of this statement.
  We now show that every graph~\(G\in\U\) has some element of~\(M\)
  as a topological minor.
  Throughout the remainder of the proof we suppose that~\(G\in\U\).
  
  If~\(G\) is connected and has three or more independent cycles, then
  the result follows from Lemma~\ref{lem: three_cycle_combinations}.
  If~\(G\) is connected and has exactly two independent cycles, then
  by Lemma~\ref{lem: two_cycle_combinations} it follows that either~\(G\)
  has an element of~\(M\) as a topological minor or~\(G\) has~\(K_4-e\) as a
  topological minor.
  For the latter case, we note that a graph with exactly two independent
  cycles having~\(K_4-e\) as a topological minor must in fact be
  a graph with two cycles sharing one or more edges along with zero
  or more leaf edges (edges incident on vertices of degree one).
  Such a graph can be reduced via \mbox{edge-contractions} to~\(K_4-e\).
  Using Theorem~\ref{thm: edge_contraction_U}, we conclude
  that~\(K_4-e\in\U\).
  However, this contradicts the result proved in
  Lemma~\ref{lem: K4-e_is_sat}.
  We therefore conclude that any connected graph~\(G\)
  having exactly two independent cycles must have some element of~\(M\)
  as a topological minor.
  
  If~\(G\) is connected and has exactly one cycle, then~\(G\) can be
  reduced to~\(C_3\) by \mbox{edge-contractions}.
  We know from Lemma~\ref{lem: C3_is_sat} that~\(C_3\notin\U\).
  Thus, Theorem~\ref{thm: edge_contraction_U} implies~\(G\notin\U\),
  which is a contradiction.
  Thus~\(G\) cannot have exactly one cycle.
  If~\(G\) is connected and has no cycles, then~\(G\) is a tree graph.
  We know from Lemma~\ref{lem: trees_are_sat} that tree graphs are not
  in~\(\U\), which is a contradiction.
  Thus~\(G\) cannot have zero cycles.
  
  Finally, suppose~\(G\in\U\) and~\(G\) is not connected.
  Then, some connected component~\(G'\) of~\(G\) is in~\(\U\).
  By our previous argument, we know that~\(G'\) has some element of~\(M\)
  as a topological minor.
  This element of~\(M\) must therefore also be a topological minor of~\(G\).
  We have showed that any graph in \(\U\), connected or not, has an element
  of~\(M\) as a topological minor.\qed\end{pf}

\begin{cor}
  A simple graph~\(G\) can support an unsatisfiable~\mbox{\(2\)-CNF} if and
  only if one of the four graphs in Lemma~\ref{lem: three_cycle_combinations}
  is a topological minor of~\(G\).  
\end{cor}
\begin{pf}
  This is simply an alternative way of stating the result of
  Theorem~\ref{thm:main} using the full definition of~\(\U\) and of the
  complete set of minimal unsatisfiability graphs.\qed
\end{pf}

\begin{rmk} As embedded graphs in \(\mathbb{R}^3\), the four graphs in
  Theorem~\ref{thm:main} are homeomorphic to one of the
  following~\mbox{\(1\)-dimensional} cell complexes

  {\centering
  \vspace*{1ex}
  {\setlength{\tabcolsep}{1em}
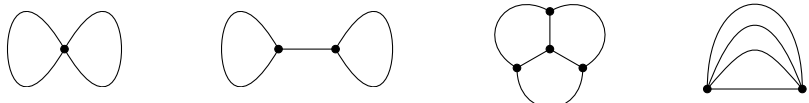
    \begin{tabular}{cccc}
      \FlexibleButterfly{}     & \FlexibleBowtie{}        & \FlexibleKfour{}             & \FlexibleBook{}\vspace*{-2ex}
  \end{tabular}}
  \vspace*{2ex}}
  \\One can immediately identify unsatisfiable sentences supported on these
  cell complexes.
  For the first one, this is the
  (unreduced)~\mbox{\(2\)-CNF}~\(({x\vee x)\wedge(\neg x\vee \neg x)}\).
  For the second cell complex, the
  (unreduced)~\(2\)-CNF \({(x\vee x)\wedge(\neg x\vee\neg y)\wedge(y\vee y)}\)
  is unsatisfiable.
  For the third, the cell complex is the same as~\(K_4\) and any
  unsatisfiable~\(2\)-CNF supported on~\(K_4\) will suffice.
  For the fourth cell complex, the (unreduced)~\mbox{\(2\)-CNF}
  \({(x\vee y)\wedge(\neg x\vee y)\wedge(x\vee\neg y)\wedge(\neg x\vee\neg y)}\)
  is unsatisfiable.
  
\end{rmk}

\section{Conclusion -- relation between our result and Robertson-Seymour 
theorem}\label{sec:conclusion}

In our main result, Theorem~\ref{thm:main}, we showed that given a
simple graph, all reduced sentences supported on it will be
satisfiable if and only if four specific graphs are forbidden as
topological minors of the original graph.
Thus it is natural to ask if the forbidden minors theory of
Robertson-Seymour~\cite{RoSe1983}--\cite{RoSe2004} applies to our
problem.

\subsection{Graph minors, pseudo-minors and topological minors}
A graph \(G\) is traditionally said to be a \emph{minor} of \(H\)
if it can be obtained from \(H\) by a series of \mbox{edge-deletions},
\mbox{vertex-deletions} and/or \mbox{edge-contractions}.  The
edge-contraction operation involves picking an edge in the graph,
deleting it, then merging the two end-vertices of that edge. The
problem is that this merging can potentially create multi-edges
(meaning, even if \(H\) is a simple graph, its minor \(G\) might be a
multigraph). There are three ways to deal with this creation of
multigraphs in our discussion:
\begin{enumerate}[a)]
\item We can \emph{allow} it: in that case, we cannot restrict
  ourselves to only dealing with simple graphs and we will have to
  consider sentences supported on multigraphs too. This can be done
  but is complicated and we consider it to be beyond the scope of this
  paper.
\item We can \emph{rectify} it: in this case, we can let the
  multi-edge be created, to only then ``throw out'' the extra
  edge. However, the trouble here is this ``throwing out''
  operation. Since we have for each \(2\)-CNF an associated
  multigraph, throwing out this extra edge corresponds to throwing out
  a clause from the~\mbox{\(2\)-CNF}. There is no real justification for
  doing this. Seen another way, throwing out an extra edge from the
  graph without justification will break the correspondence that the
  graph has with the \(2\)-CNF it supports.
\item We can \emph{forbid} it: and this is the path we choose in this
  paper. To forbid the creation of multigraphs we allow
  edge-contractions but only at edges not contained in triangles.
\end{enumerate}

Hence the relation we consider is not that of a ``minor'' but that of
a ``pseudo-minor'' (edge-contractions are allowed but only at edges
not contained in triangles).  We note here that pseudo-minor is not a
standard term used among graph theorists.

\subsection{Robertson-Seymour graph minor theorem}
The Robertson-Seymour graph Theorem (RST) states that if a
\mbox{graph-family} \(\mathcal{F}\) is \mbox{minor-closed}, meaning
that every minor of a graph in~\(\mathcal{F}\) is also
in~\(\mathcal{F}\), then there are at most
finitely-many forbidden minors of~\(\mathcal{F}^c\) -- meaning a graph
belongs to~\(\mathcal{F}\) if and only if it does not contain as a
minor any of these forbidden graphs.
In our case, a suitable candidate for~\(\mathcal{F}\) is the
set of graphs that support only satisfiable~\(2\)-CNFs.
According to our notation, this is \(\U^c\)).  If~\(\U^c\) is shown to be
minor-closed then by RST we will be able to find a finite number of
forbidden minors.  However, as discussed above, we cannot meaningfully
talk about minors of graphs in the context of~\(2\)-CNFs while
restricting ourselves to simple graphs only.  We deal with the
relation of a ``pseudo-minor'' and hence RST does not apply.
Unfortunately, the authors of this paper are not aware of a
Robertson-Seymour-type theorem for pseudo-minors.

To summarize, if~\(\U^c\) were shown to be minor-closed, then RST
would imply finitely-many forbidden minors.  In this paper we have
established that~\(\U^c\) is pseudo-minor-closed.  However, RST does
not apply so it appears as though it is possible that there might be
infinitely-many forbidden pseudo-minors.
Fortunately, our result implies that are in fact only finitely many.

We have also shown that the graph family~\(\U^c\) is
\mbox{topological-minor-closed}.
However, the topological minor relation is not a \mbox{well-quasi-ordering}
on the set of finite graphs and therefore, RST does not apply to it.
Fortunately, our main theorem established
that~\(\U^c\) has exactly four forbidden topological minors
(we call this forbidden set the complete set of minimal unsatisfiability graphs).

\subsection{Computational complexity}
There are subtle results involving the computational complexity
implications of our result as well.  Since~\mbox{\(2\)-SAT} is known
to be in complexity class P, it would at first glance appear that we
have made the problem worse: one can infer from our main theorem that the
satisfiability of a given \(2\)-CNF can be established by looking for
four specific topological minors in its associated graphs.
If these four topological minors are absent, then the \(2\)-CNF must
be satisfiable.
If they are present, then the \(2\)-CNF may or may not be satisfiable.
This may appear to be setback since the \textsc{graph minor} decision
problem is known to be NP-complete, implying that the \textsc{graph
  topological minor} problem will also be NP-complete.

However, on closer inspection, one realizes that the decision problem
of determining if a fixed graph is present as a topological minor can
actually be resolved in polynomial time as a consequence of the graph
minor theorem.
Since our set of forbidden topological minors is finite we can therefore
guarantee that the task of searching for them can be accomplished in
polynomial time.

Thus, we have done ``no worse'' than the original \(2\)-SAT problem as
both can still be solved in polynomial time. We emphasize however that
as we have stated earlier, our goal in this paper is the
characterization of unsatisfiable reduced \(2\)-CNFs as opposed to
questions of algorithmic efficiency for solving the satisfiability
problem for \(2\)-CNFs for which, in any case, efficient algorithms
already exist. We hope that the techniques developed in this paper
when generalized to hypergraphs might shed light on the structure of
unsatisfiable~\mbox{\(3\)-CNFs}.

\paragraph{Acknowledgement:} The authors thank Yuliy Baryshnikov for
suggesting the problem of studying satisfiability from the viewpoint
of the underlying structure of the sentences and for early discussions
on the subject.

\bibliographystyle{elsarticle-num.bst}
\bibliography{sat.bib}

\end{document}